\documentclass[a4paper,12pt]{article}     

\usepackage{graphicx}		  
\usepackage{amsmath,amssymb}	

\usepackage{makeidx}

\newtheorem{lem}{Lemma}[section]
\newtheorem{thm}[lem]{Theorem}
\newtheorem{proposition}[lem]{Proposition}

\newtheorem{assumption}{Assumption}
\newtheorem{remark}[lem]{Remark}

\def\authorfont{\footnotesize}

\def\ccode#1{\par		
	\vspace*{8pt}
	{\authorfont{\leftskip18pt\rightskip\leftskip
	\noindent #1\par}}\par}

\newenvironment{Proof}{
\hspace*{-9mm}
{ \it Proof.}}
{\hfill {$\square$}\vspace{1.5em}}

\begin{document}

\begin{center}{
{\Large The structure of a minimal $n$-chart with 
two crossings I:\\
Complementary domains of $\Gamma_1\cup\Gamma_{n-1}$
}
\vspace{10pt}
\\ 
Teruo NAGASE and Akiko SHIMA\footnote{The second author was supported by JSPS KAKENHI Grant Number 18K03309.}\\
Department of Mathematics, 
Tokai University,\\
4-1-1 Kitakaname, Hiratuka, Kanagawa, 259-1292 Japan.\\
email: nagase@keyaki.cc.u-tokai.ac.jp,
shima@keyaki.cc.u-tokai.ac.jp
}

\end{center}

\begin{abstract}
This is the first step of the two steps 
to enumerate 
the minimal charts with two crossings. 
For a label $m$ of a chart $\Gamma$
we denote by $\Gamma_m$ 
the union of all the edges 
of label $m$ and their vertices. 
For a minimal chart $\Gamma$ 
with exactly two crossings, 
we can show that 
the two crossings are contained in 
$\Gamma_\alpha\cap\Gamma_\beta$ 
for some labels $\alpha<\beta$.
In this paper, we 
study the structure of 
a disk $D$ 
not containing any crossing 
but satisfying 
$\Gamma\cap \partial D
\subset 
\Gamma_{\alpha+1}\cup \Gamma_{\beta-1}$.
\end{abstract}

\ccode{2010 Mathematics Subject Classification. Primary 57Q45; Secondary 57Q35.}
\ccode{ {\it Key Words and Phrases}. surface link, chart, crossing. }

\setcounter{section}{0}
\section{{\large Introduction}}


Charts are oriented labeled graphs
in a disk 
with three kinds of vertices
called black vertices, crossings,
and white vertices (see Section~\ref{s:Prel} for the precise definition of charts, 
black vertices, crossings, and white vertices).
From a chart, we can construct an oriented closed surface 
embedded in 4-space ${\Bbb R}^4$ 
 (see \cite[chapter 14, chapter 18 and chapter 23]{BraidBook}). 
A C-move 
is a local modification between two charts
in a disk (see Section~\ref{s:Prel}).
A C-move between two charts induces 
an ambient isotopy between oriented closed surfaces 
corresponding to the two charts.
Two charts are said to be {\it C-move equivalent}\index{C-move~equivalent} 
if there exists
a finite sequence of C-moves 
which modifies one of the two charts 
to the other.

We will work in the PL or smooth category. 
All submanifolds are assumed to be locally flat.
A {\it surface link} is a closed surface embedded in 4-space ${\Bbb R}^4$. 
A {\it $2$-link} is a surface link each of whose connected component is a $2$-sphere.
A {\it $2$-knot}
is a surface link which is a $2$-sphere.
An orientable surface link is called a 
{\it ribbon surface link}
if there exists an immersion of a 3-manifold $M$
into ${\Bbb R}^4$ sending the boundary of $M$ onto the surface link
such that each connected component of $M$ is a handlebody
and its singularity
consists of ribbon singularities,
here a ribbon singularity
is a disk in the image of $M$
whose pre-image consists of 
two disks;
one of the two disks is a proper disk of $M$ 
and
the other is a disk in the interior of $M$.
In the words of charts,
a ribbon surface link is
a surface link corresponding to a {\it ribbon chart}, 
a chart C-move equivalent to 
a chart
without white vertices \cite{BraidThree}.
A chart is called a {\it $2$-link chart}
if a surface link corresponding to the chart is a $2$-link.

In this paper, 
we denote the closure, the interior, 
the boundary, and the complement of $(...)$ by 
$Cl(...)$, Int$(...)$, 
$\partial(...)$, $(...)^c$ 
respectively.
Also for a finite set $X$, 
the notation $|X|$ denotes 
the number of elements in $X$.

At the end of this paper 
there are lists of terminologies and notations 
which are used in this paper.

Kamada showed that 
any $3$-chart is a ribbon chart 
\cite{BraidThree}.
Kamada's result was extended by Nagase and Hirota:
Any $4$-chart with at most one crossing
is a ribbon chart \cite{NH}.
We showed that any $n$-chart with at most one crossing is a ribbon chart
\cite{OneCrossing}.
We also showed that any $2$-link chart 
with at most two crossings
 is a ribbon chart  \cite{TwoCrossingI},  
 \cite{TwoCrossingII}.

\index{$\Gamma_m$}
Let $\Gamma$ be a chart.
For each label $m$, we define
$$\Gamma_m=\text{ 
the union of 
all the edges of label $m$ and 
their vertices in }\Gamma.$$

Let $\Gamma$ be a chart in a disk $D^2$, and 
$D$ a disk in $D^2$.
The pair $(\Gamma\cap D,D)$ is called 
a {\it tangle}\index{tangle} provided that 
\begin{enumerate}
\item[(i)]
$\partial D$ does not contain
any white vertices, 
black vertices
nor crossings of $\Gamma$, 
\item[(ii)]
if an edge of $\Gamma$ intersects $\partial D$, 
then the edge intersects $\partial D$ transversely, 
\item[(iii)] $\Gamma\cap D\not=\emptyset$.
\end{enumerate}

In this paper and \cite{StII},
we investigate 
the structure of minimal charts with two crossings
(see Section~\ref{s:Prel}
for the precise definition of a minimal chart),
and give an enumeration of the charts 
with two crossings. 
First, we split a minimal chart with two crossings 
into two kinds of tangles; 
one is called a net-tangle, and 
the other is called an IO-tangle. 

\begin{figure}
\begin{center}
\includegraphics{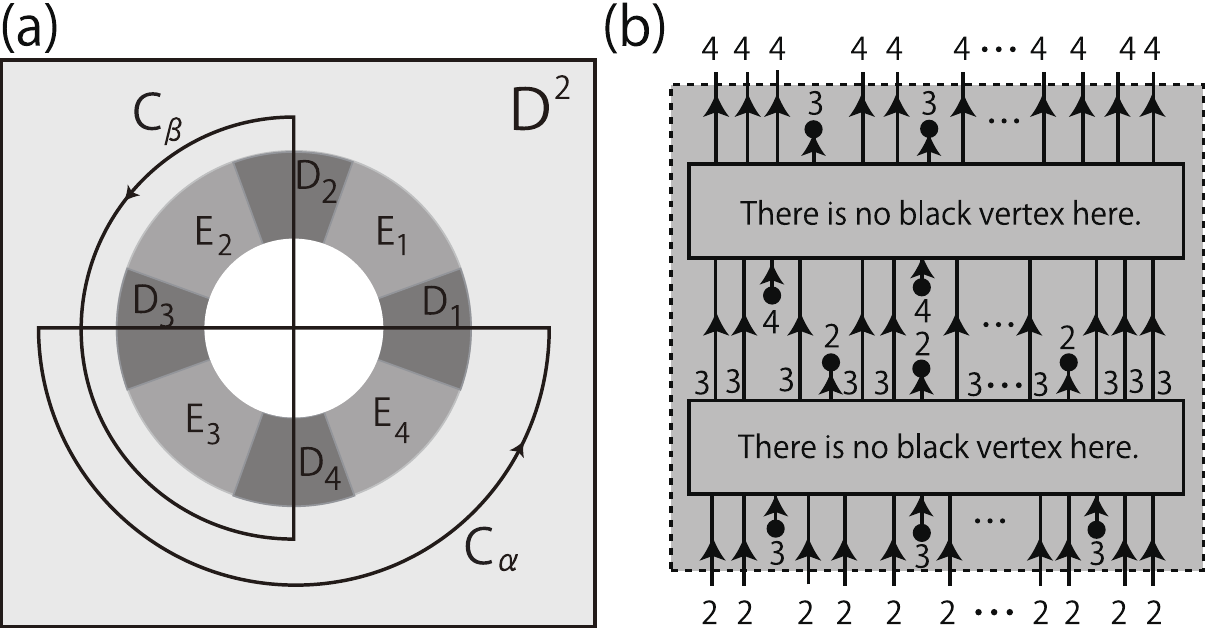}
\end{center}
\caption{\label{fig-01}
(b) a tangle $(\Gamma\cap E_i,E_i)$ with $\Gamma\cap E_i\subset \Gamma_2\cup \Gamma_3\cup \Gamma_4$
for the case $\alpha=1$ and $\beta=5$,
here all the free edges and simple hoops are in a regular neighbourhood of $\partial D^2$,
and the numbers are labels of the chart.
}
\end{figure}

\begin{figure}[htb]
\begin{center}
\includegraphics{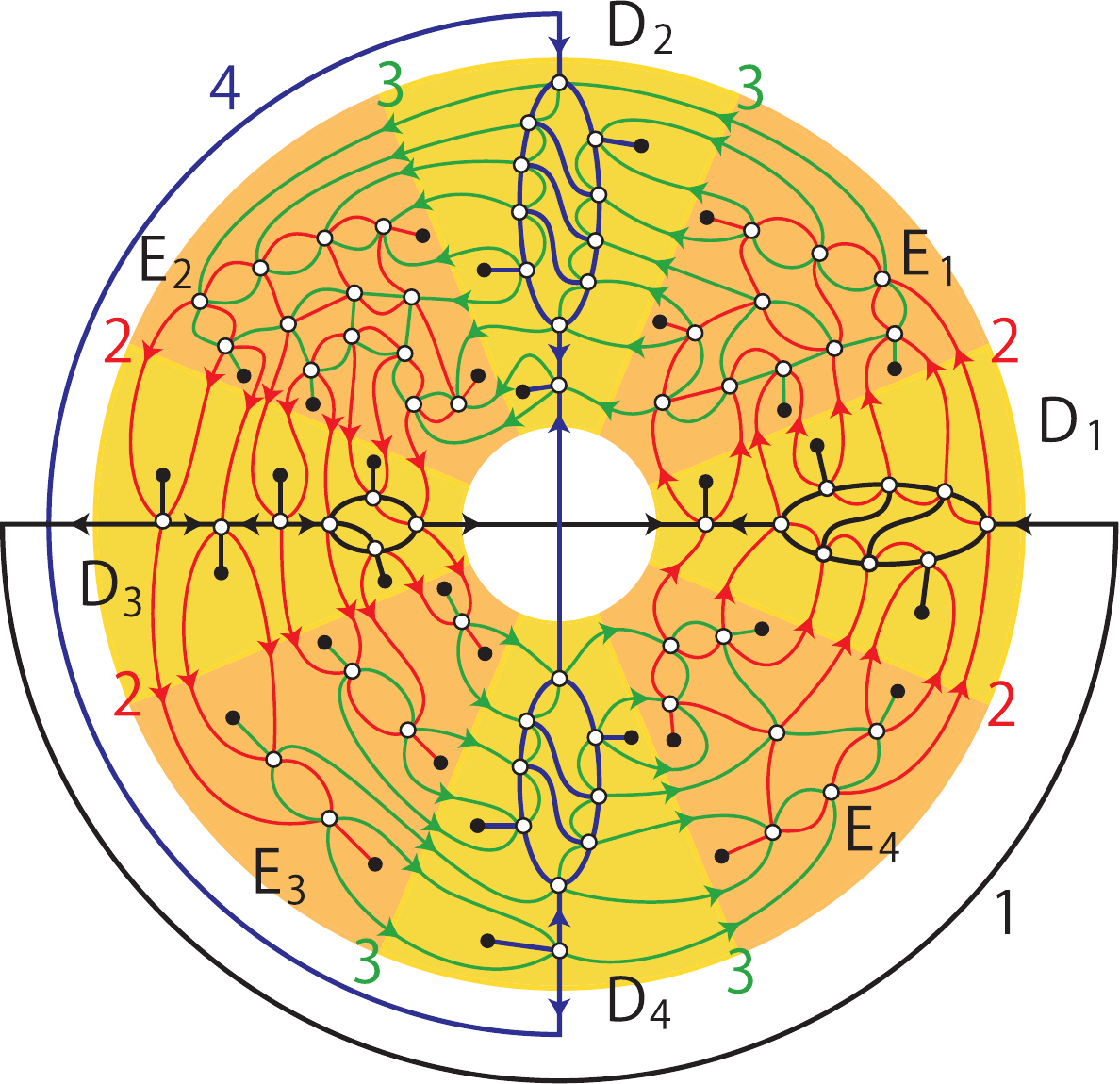}
\end{center}
\caption{ \label{fig-02} A 5-chart with two crossings.}
\end{figure}

We investigate net-tangles in this paper, 
and
IO-tangles in \cite{StII}.
In short, for any minimal $n$-chart $\Gamma$ 
with two crossings in a disk $D^2$,
there exist two labels 
$1\le \alpha<\beta\le n-1$ such that 
 $\Gamma_\alpha$ and $\Gamma_\beta$
contain cycles $C_\alpha$ and $C_\beta$
with $C_\alpha\cap C_\beta$ the two crossings
and that
for any label $k$ with $k<\alpha$ or $\beta<k$,
the set $\Gamma_k$ does not contain a white vertex.
If $\Gamma_\alpha$ or $\Gamma_\beta$
contains at least three white vertices,
then
after shifting all the free edges and simple hoops (see Section~\ref{s:Prel} for the definition of free edges, hoops, and simple hoops)
into a regular neighbourhood of $\partial D^2$
by applying C-I-M1 moves and C-I-M2 moves, 
we can find 
an annulus $A$ 
containing all the white vertices of $\Gamma$ 
but not intersecting 
any hoops nor free edges
such that (see Fig.~\ref{fig-01}(a))
\begin{enumerate}
\item[(1)] 
each connected component of $Cl(D^2-A)$ 
contains a crossing,
\item[(2)]
$\Gamma\cap \partial A=
(C_\alpha\cup C_\beta)\cap \partial A$, and 
$\Gamma\cap \partial A$ 
consists of eight points.
\end{enumerate}
We can show 
the annulus $A$ can be split 
into 
mutually disjoint four disks 
$D_1,D_2,D_3,D_4$ and 
mutually disjoint four disks 
$E_1,E_2,E_3,E_4$ such that
\begin{enumerate}
\item[(3)] 
for each $i=1,2,3,4$, 
if $\alpha+1=\beta-1$, 
then $\Gamma\cap E_i$ consists of 
parallel proper arcs of $E_i$ of label $\alpha+1$, 
otherwise 
the tangle $(\Gamma\cap E_i,E_i)$ is 
a net-tangle with 
$\Gamma\cap E_i\subset
\cup^{\beta-1}_{j=\alpha+1}\Gamma_j$
as shown in Fig.~\ref{fig-01}(b), and
\item[(4)] 
for each $i=1,3$ (resp. $i=2,4$)
the tangle $(\Gamma\cap D_i,D_i)$ 
is an IO-tangle of label $\alpha$ 
(resp. label $\beta$).
\end{enumerate}
We count the number of edges between edges 
with black vertices
in Fig.~\ref{fig-01}(b) 
to enumerate charts with two crossings.
As important results, 
from the enumeration 
we can calculate the fundamental group of 
the exterior of the surface link represented by $\Gamma$, 
and the braid monodromy of the surface braid 
represented by $\Gamma$.
For example, 
the normal form
for the 5-chart in Fig.~\ref{fig-02} (see \cite{StII})
is\\ 
$((5,3,1,2),(8,9,8,7);\\
(1,3,4),(3,3,2),(2,3,3),(6,2);\\
(1,3,5),(2,3,4),(4,3,2),(2,3,3,1);\\
(2,4,2),(1,3,2,2),(1,2,4,1),(4,3,1);\\
(2,5),(2,3,2),(2,3,2),(5,2))$.

If an edge $e$ of a chart $\Gamma$ is oriented 
from a vertex $v_1$ to the other vertex $v_2$, 
then we say that 
the edge $e$ is {\it outward}\index{outward} at $v_1$, and 
the edge $e$ is {\it inward}\index{inward} at $v_2$.

\begin{figure}
\begin{center}
\includegraphics{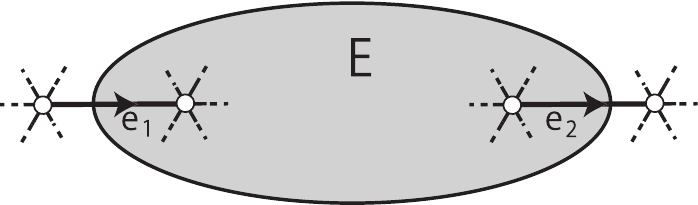}
\end{center}
\caption{ \label{fig-03}
The gray area is a disk $E$.
The edge $e_1$ is an I-edge for $E$, and
the edge $e_2$ is an O-edge for $E$.
}
\end{figure}

Let $\Gamma$ be a chart, and 
$E$ a disk.
An edge $e$ of the chart $\Gamma$
is called 
an {\it I-edge $($resp. O-edge$)$}\index{I-edge} \index{O-edge}
for $E$ provided that (see Fig.~\ref{fig-03})
\begin{enumerate}
\item[(i)] 
the edge $e$ possesses two white vertices,
one is in Int~$E$  
and 
the other in $E^c$,
\item[(ii)] 
the edge $e$ intersects $\partial E$ 
by exactly one point, and
\item[(iii)] 
the edge $e$ is inward (resp. outward) at
the vertex in Int~$E$.
\end{enumerate}
We often say just {\it an I-edge} 
instead of {\it an I-edge for $E$}
if there is no confusion.
Similarly we often say just {\it an O-edge} 
instead of {\it an O-edge for $E$}.

For a simple arc $\ell$, we set\\
\ \ \ $
\begin{array}{rl}
\partial \ell&= \text{ the set of its two endpoints, 
and}\\
{\rm Int}~\ell&=~\ell-\partial \ell.
\end{array}
$

\begin{figure}
\begin{center}
\includegraphics{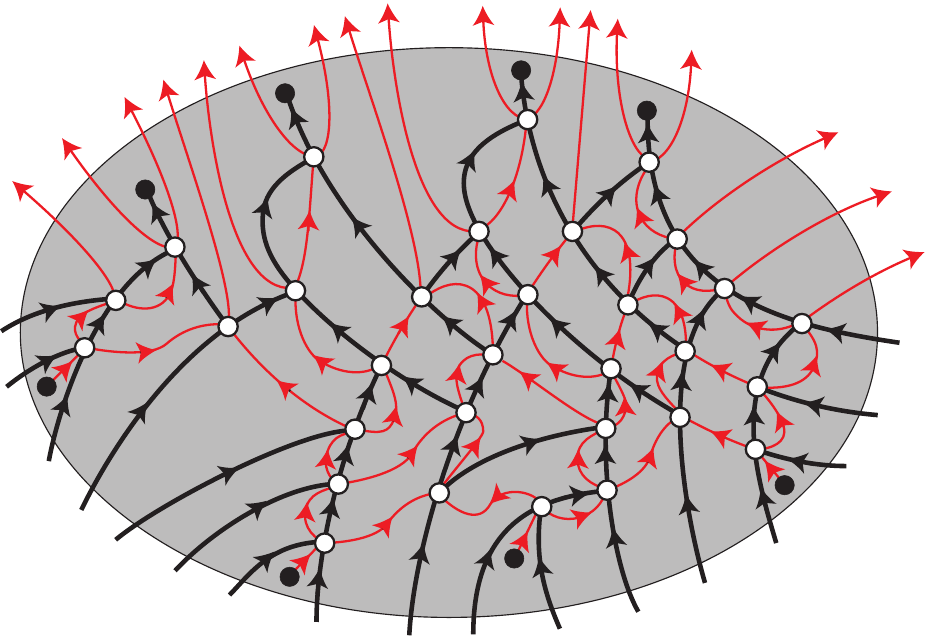}
\end{center}
\caption{ \label{fig-04}
An example of an upward net-tangle with 
a label pair $(\alpha,\alpha+1)$.
The gray area is a disk $D$,
the thick edges are of label $\alpha$, and
the thin edges are of label $\alpha+1$.
}
\end{figure}

Let $\Gamma$ be a chart.
A tangle $(\Gamma\cap D,D)$ 
is called a {\it net-tangle}\index{net-tangle}
provided that
\begin{enumerate}
\item[(i)]
the disk $D$ contains no crossing, hoop, 
nor free edge
but a white vertex,
\item[(ii)]
there exist two labels $\alpha,\beta$ 
with $\alpha<\beta$ and
$\Gamma\cap D
\subset
\cup_{i=\alpha}^{\beta}\Gamma_{i}$,
and
\item[(iii)]
there exist two arcs 
$L_\alpha,L_\beta$ on $\partial D$ 
with $L_\alpha\cap L_\beta$ two points such that
\begin{enumerate}
\item[(a)] 
$\Gamma\cap\partial D=
(\Gamma_\alpha\cap {\rm Int}~L_\alpha)
\cup 
(\Gamma_\beta\cap {\rm Int}~L_\beta)$,
\item[(b)]
all the edges intersecting $L_\alpha$ 
are I-edges of label $\alpha$ or\\
all the edges intersecting $L_\alpha$ 
are O-edges of label $\alpha$,
and
\item[(c)]
all the edges intersecting $L_\beta$ 
are O-edges of label $\beta$ or\\
all the edges intersecting $L_\beta$ 
are I-edges of label $\beta$.
\end{enumerate}
\end{enumerate}
The pairs $(\alpha,\beta)$ and 
$(L_\alpha,L_\beta)$ are called 
a {\it label pair}\index{label pair} 
and 
a {\it boundary arc pair}\index{boundary arc pair} of the net-tangle
respectively.
If all the edges of label $\alpha$ 
intersecting $L_\alpha$
are I-edges (resp. O-edges), 
and
if all the edges of label $\beta$ 
intersecting $L_\beta$
are 
O-edges (resp. I-edges),
then the net-tangle is said to be 
{\it upward} (resp. {\it downward}) \index{upward net-tangle}\index{downward net-tangle}
(see Fig.~\ref{fig-04}).
An upward or downward 
net-tangle with 
a label pair $(\alpha,\alpha+1)$ 
is called an N-{\it tangle}.\index{N-tangle}

An edge of a chart $\Gamma$ is called a {\it terminal edge}\index{terminal edge}
if it contains a white vertex and a black vertex.
If 
a terminal edge is inward
at its black vertex,
then the edge 
is called an {\it I-terminal edge}\index{I-terminal edge},
otherwise the edge is called an {\it O-terminal edge}.\index{O-terminal edge} 

For minimal charts,
we shall show the following two theorems:

\begin{thm}\label{Theorem1}
Let $\Gamma$ be a minimal chart, and 
$(\Gamma\cap D,D)$ a net-tangle
with a label pair $(\alpha,\alpha+1)$. 
Then we have the following.
\begin{enumerate}
\item[{\rm (a)}] 
The tangle is an N-tangle.
\item[{\rm (b)}] 
The number of the edges of label $\alpha$ 
intersecting $\partial D$
is equal to
the number of the edges of label $\alpha+1$ 
intersecting $\partial D$.
\item[{\rm (c)}] 
There exists a terminal edge in $D$.
\item[{\rm (d)}]
The number of terminal edges of 
label $\alpha$ in $D$ 
is equal to 
the number of terminal edges of 
label $\alpha+1$ in $D$.
\item[{\rm (e)}]
If the tangle is upward $($resp. downward$)$, 
then 
all the terminal edges of label $\alpha$ in $D$ 
are I-terminal $($resp. O-terminal$)$ edges
and 
all the terminal edges of label $\alpha+1$ in $D$ 
are O-terminal $($resp. I-terminal$)$ edges.
\end{enumerate}
\end{thm}

Let $E$ be a connected compact surface.
A simple arc $\ell$ in $E$ is called a 
{\it proper arc}\index{proper arc} provided that 
$\ell\cap \partial E=\partial \ell$.

There exists a special C-move 
called a C-I-M2 move 
(see Fig.~\ref{fig-06} in Section~\ref{s:Prel} 
for C-I-M2 moves).
Let $\Gamma$ be a minimal chart, and 
$(\Gamma\cap D,D)$ 
a net-tangle with 
a label pair $(\alpha,\beta)$. 
Let $\widetilde\Gamma$ be a minimal chart 
such that $(\widetilde\Gamma\cap D, D)$ 
is a tangle without crossing nor hoops.
Then the chart $\widetilde\Gamma$ is said to be 
M2-{\it related to $\Gamma$ with respect to $D$}\index{M2-related to a chart}
provided that
\begin{enumerate}
\item[(i)]
$\Gamma\cap D^c=\widetilde\Gamma\cap D^c$, and
\item[(ii)]
the chart $\widetilde\Gamma$ is obtained from
the chart $\Gamma$ by a finite sequence of 
C-I-M2 moves in $D$ each of which modifies 
two edges of some label $i$ 
with $\alpha<i<\beta$.
\end{enumerate}

\begin{thm}\label{Theorem2}
Let $\Gamma$ be a minimal chart, and 
$(\Gamma\cap D,D)$ a net-tangle
with a label pair $(\alpha,\beta)$. 
Then $\Gamma$ is M2-related to 
a minimal chart $\widetilde\Gamma$ 
with respect to $D$ such that
there exist N-tangles 
$(\widetilde\Gamma\cap D_{\alpha},D_{\alpha}),
(\widetilde\Gamma\cap D_{\alpha+1},D_{\alpha+1}),
\cdots,
(\widetilde\Gamma\cap D_{\beta-1},D_{\beta-1})$
equipped with
\begin{enumerate}
\item[{\rm (a)}] for each 
$i=\alpha,\alpha+1,\cdots,\beta-1$,
the tangle 
$(\widetilde\Gamma\cap D_{i},D_{i})$ 
is an N-tangle with the label pair $(i,i+1)$,
\item[{\rm (b)}] $D=\cup_{i=\alpha}^{\beta-1}D_i$, 
\item[{\rm (c)}] 
for each $i=\alpha,\alpha+1,\cdots,\beta-2$,
the intersection $D_i\cap D_{i+1}$ 
is a proper arc of $D$,
\item[{\rm (d)}]
all the N-tangles are 
upward or downward
simultaneously.
\end{enumerate}
\end{thm}

Our paper is organized as follows:
In Section~\ref{s:Prel}, 
we introduce 
the definition of charts and its related words.
In Section~\ref{s:OneWayPath},
we investigate a directed path of label $m$.
In Section~\ref{s:Ntangle},
we investigate terminal edges and directed paths of label $m$ in an N-tangle.
In Section~\ref{s:Spindle},
we investigate  
a disk $E$ such that $\partial E$ consists of two directed paths of label $m$.
In Section~\ref{s:HalfSpindle},
we investigate  
a disk $E$ such that $\partial E$ consists of an arc and two directed paths of label $m$.
In Section~\ref{s:ProofTh1},
we prove Theorem~\ref{Theorem1}.
In Section~\ref{s:ProofTh2},
we prove Theorem~\ref{Theorem2}.

\section{{\large Preliminaries}}
\label{s:Prel}

In this section, 
we introduce 
the definition of charts and its related words.

\begin{figure}[htb]
\begin{center}
\includegraphics{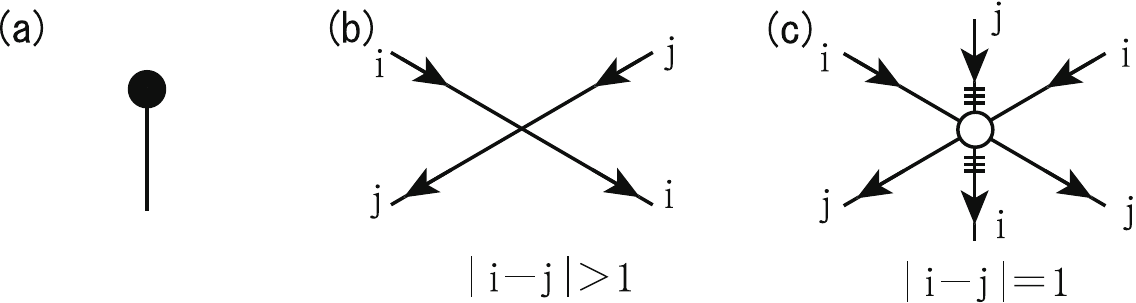}
\end{center}
\caption{ \label{fig-05} (a) a black vertex. 
(b) a crossing. (c) a white vertex. 
Each arc with three transversal short arcs is a middle arc at the white vertex. }
\end{figure}

Let $n$ be a positive integer.
An $n$-{\it chart}  
(a braid chart of degree $n$ \cite{KS}
or a surface braid chart of degree $n$ \cite{BraidBook}) 
is 
an oriented labeled graph in the interior of a disk,
which may be empty 
or
have closed edges without vertices
satisfying the following four conditions
(see Fig.~\ref{fig-05}):
\begin{enumerate}
\item[(i)] 
every vertex has degree $1$, $4$, or $6$.
\item[(ii)] 
the labels of edges are 
in $\{1,2,\dots,n-1\}$.
\item[(iii)]
in a small neighborhood of
each vertex of degree $6$,
there are six short arcs,
three consecutive arcs are
oriented inward 
and
the other three are outward,
and
these six are labeled $i$ and $i+1$
alternately for some $i$,
where the orientation and label of
each arc are inherited from
the edge containing the arc.
\item[(iv)]
for each vertex of degree $4$,
diagonal edges have the same label
and
are oriented coherently,
and the labels $i$ and $j$ of
the diagonals satisfy $|i-j|>1$.
\end{enumerate}
We call a vertex of degree $1$ a {\it black vertex},
a vertex of degree $4$ a {\it crossing}, and 
a vertex of degree $6$ a {\it white vertex}
respectively.
Among six short arcs
in a small neighborhood of
a white vertex,
a central arc of each three consecutive arcs
oriented inward (resp. outward) 
is called a   
{\it middle arc}\index{middle arc} at the white vertex
(see Fig.~\ref{fig-05}(c)).
For each white vertex $v$, 
there are two middle arcs at $v$ 
in a small neighborhood of $v$.
An edge $e$ is said to be 
{\it middle at}\index{middle at $v$} a white vertex $v$ if it contains a middle arc at $v$.


\begin{figure}[htb]
\begin{center}
\includegraphics{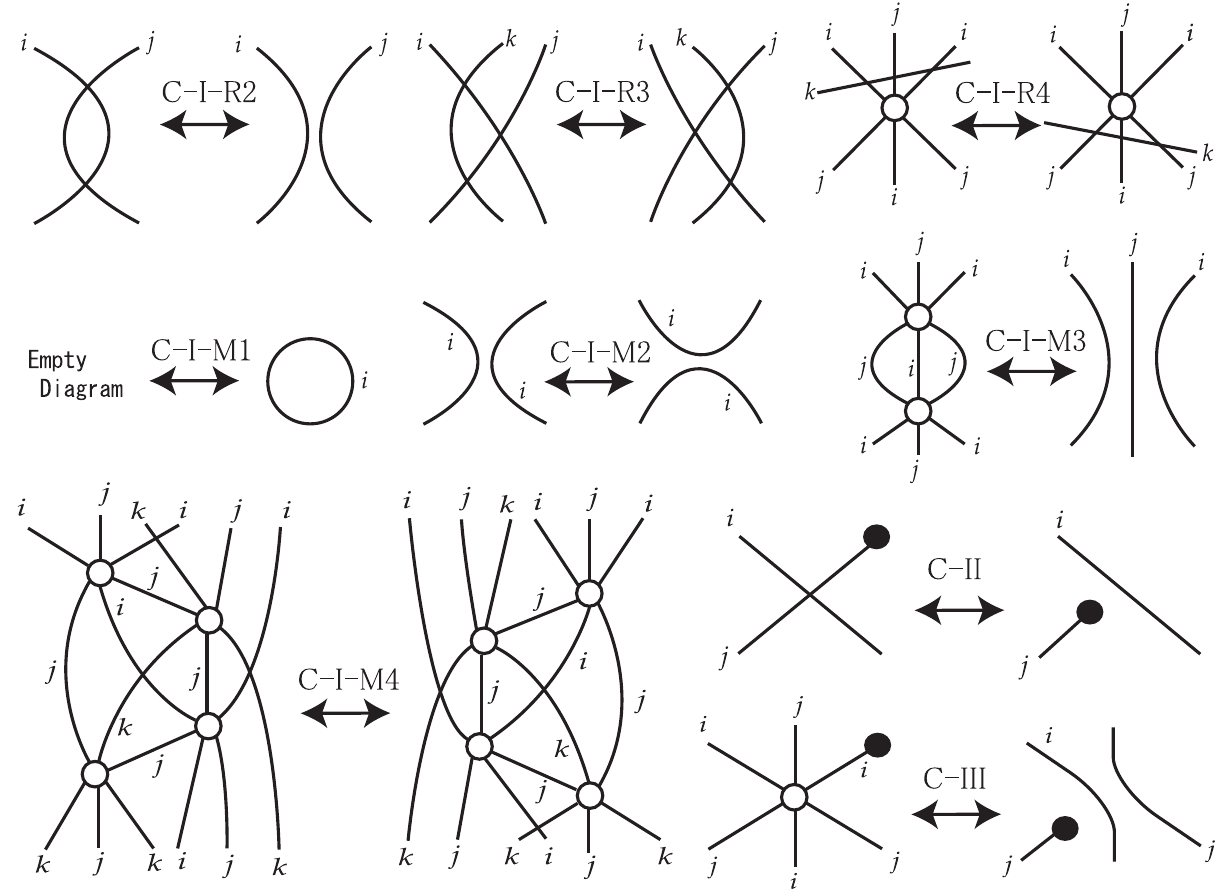}
\end{center}
\caption{ \label{fig-06} For the C-III move, the edge containing the black vertex does not contain a middle arc at
a white vertex in the left figure. }
\end{figure}

Now {\it C-moves} are local modifications 
of charts as shown in Fig.~\ref{fig-06}
(cf. \cite{KS}, \cite{Braid},
\cite{BraidBook} and \cite{Tanaka}).
We often use C-I-M2 moves, C-I-M3 moves, C-II moves
and C-III moves. 

Kamada originally defined CI-moves
as follows: 
A chart $\Gamma$ is obtained from
a chart $\Gamma'$ in a disk $D^2$
by a {\it CI-move},
if there exists a disk $E$ 
in $D^2$ such that 
\begin{enumerate}
\item[(i)] 
the two charts $\Gamma$ and $\Gamma'$
intersect the boundary of $E$ transversely
or
do not intersect the boundary of $E$, 
\item[(ii)] 
$\Gamma\cap E^c=\Gamma'\cap E^c$, and
\item[(iii)]
neither $\Gamma\cap E$ nor 
$\Gamma'\cap E$ 
contains a black vertex,
\end{enumerate}
where $E^c$ is 
the complement of $E$ in the disk $D^2$.

\begin{remark}
{\rm Any CI-move is realized by a finite sequence of 
seven types: C-I-R2, C-I-R3, C-I-R4, 
C-I-M1, C-I-M2, C-I-M3, C-I-M4.}
\end{remark}

An edge of a chart is called 
a {\it free edge}\index{free edge}
if it contains
two black vertices.

\begin{figure}
\begin{center}
\includegraphics{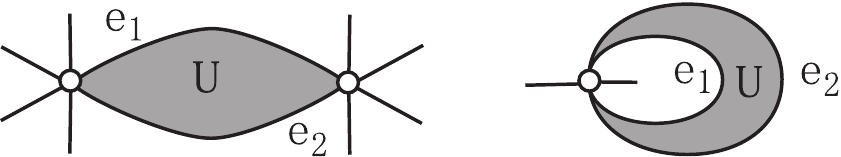}
\end{center}
\caption{ \label{fig-07} Bigons.}
\end{figure}

Let $\Gamma$ be a chart. 
Let $e_1$ and $e_2$ be edges of $\Gamma$
which connect two white vertices $w_1$ and $w_2$
where possibly $w_1=w_2$.
Suppose that 
the union $e_1\cup e_2$ bounds 
an open disk $U$.
Then $Cl(U)$ 
is called 
a {\it bigon}\index{bigon} of $\Gamma$
provided that
any edge containing $w_1$ or $w_2$ 
does not intersect the open disk $U$
(see Fig.~\ref{fig-07}).
Since $e_1$ and $e_2$ are edges of $\Gamma$,
 neither $e_1$ nor $e_2$ contains a crossing.

%

Let $\Gamma$ be a chart.
Let $w(\Gamma), f(\Gamma), 
c(\Gamma)$, and
$b(\Gamma)$ be 
the number of white vertices of $\Gamma$, 
the number of free edges of $\Gamma$,
the number of crossings of $\Gamma$,
and
the number of bigons of $\Gamma$
respectively.
The 4-tuple $(c(\Gamma),w(\Gamma),-f(\Gamma),-b(\Gamma))$ is called a 
{\it $c$-complexity} of the chart $\Gamma$.
The 4-tuple $(w(\Gamma),c(\Gamma),-f(\Gamma),-b(\Gamma))$ is called a 
{\it $w$-complexity} of the chart $\Gamma$.
The 3-tuple $(c(\Gamma)+w(\Gamma),-f(\Gamma),-b(\Gamma))$ is called a 
{\it $cw$-complexity} of the chart $\Gamma$
(see \cite{BraidThree} 
for complexities of charts).

A chart $\Gamma$ is said to be 
{\it $c$-minimal $($resp. $w$-minimal or $cw$-minimal$)$} if
its $c$-complexity (resp. $w$-complexity or $cw$-complexity) is minimal among the charts 
which are C-move equivalent to 
the chart $\Gamma$
with respect to 
the lexicographical order of the 
4-tuple (or 3-tuple) of the integers.

{\it In this paper, 
if a chart is $c$-minimal, $w$-minimal or $cw$-minimal,}\index{minimal chart}

{\it then we say that the chart is minimal.}

A {\it hoop}\index{hoop} is a closed edge of a chart $\Gamma$
that contains neither crossings nor white vertices.
Therefore a hoop decomposes $\Gamma$ into disjoint pieces:
an inside, an outside and itself.
A hoop is said to be {\it simple}\index{simple hoop} 
if one of the complementary domains
of the hoop
does not contain any white vertices.
An {\it oval nest}\index{oval nest} is a free edge together 
with some concentric simple hoops.

\begin{proposition}
\label{MoveOutFreeEdge}
Let $\Gamma$ be a chart in a disk $D^2$.
For any regular neighbourhood $N$ 
of $\partial D^2$ in $D^2$,
there exists a chart $\Gamma'$ obtained from $\Gamma$ 
by C-I-M2 moves and ambient isotopies of $D^2$ 
without changing  
the complexity such that
\begin{enumerate}
\item[{\rm (a)}]
$\Gamma'\cap (D^2-N)$ contains no free edge,
\item[{\rm (b)}]
$\Gamma'\cap N$ consists of oval nests, 
simple hoops and free edges.
\end{enumerate}
\end{proposition}

\begin{figure}
\begin{center}
\includegraphics{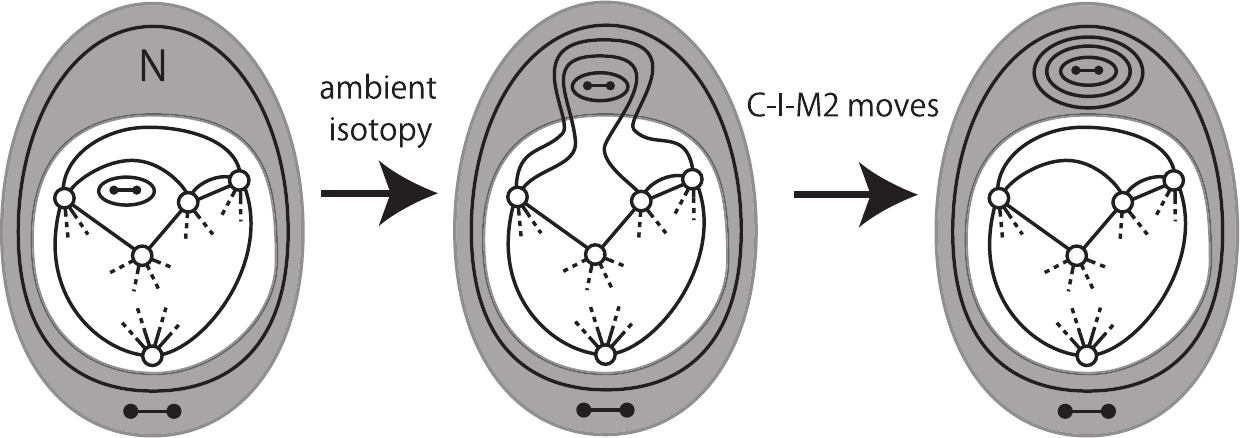}
\end{center}
\caption{ \label{fig-08} Moving free edges into 
a gray annulus $N$, a regular neighbourhood of $\partial D^2$.}
\end{figure}

\begin{Proof}
One by one, we can move free edges into 
a regular neighbourhood $N$ of $\partial D^2$ 
by C-I-M2 moves and ambient isotopies of $D^2$
as shown in Fig.~\ref{fig-08} 
satisfying (b)
(see \cite[p 238, Figure 29.2]{BraidBook}).
\end{Proof}

%


Let $m$ be a label of a chart $\Gamma$.
A simple closed curve in $\Gamma_m$
is called a {\it ring},\index{ring}
if it contains at least one crossing 
but does not contain a white vertex
nor black vertex.

Let $\Gamma$ be a chart, and 
$m$ a label of the chart. 
Let $\mathcal W$ be the set of 
all the white vertices of $\Gamma$. 
The closure of a connected component of $\Gamma_m-\mathcal W$ 
is called an {\it internal} edge of label $m$ 
if it contains a white vertex 
but does not contain any black vertex, 
here we consider $\Gamma_m$ as a topological set. 
A simple arc is said to be {\it internal} 
if it is contained in an internal edge or a ring.

\begin{proposition}
\label{Assumption0}
{\rm (\cite[remark 2.3]{MinimalChart})}
Let $\Gamma$ be a minimal chart 
in a disk $D^2$. 
Then we have 
the following:
\begin{enumerate}
\item[{\rm (a)}]
if an edge of $\Gamma$ contains a black vertex, 
then the edge is a terminal edge or a free edge.
\item[{\rm (b)}]
any terminal edge of $\Gamma$ contains a middle arc 
at its white vertex.
\item[{\rm (c)}] each complementary domain of any ring must contain at least one white vertex.
\end{enumerate}
\end{proposition}

\begin{figure}
\begin{center}
\includegraphics{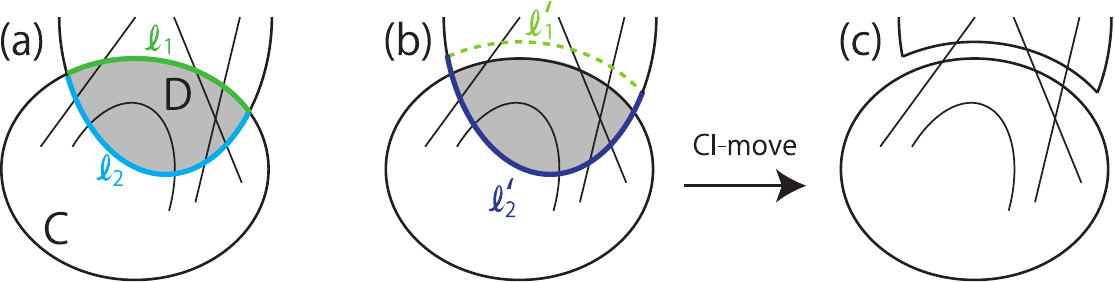}
\end{center}
\caption{ \label{fig-09} 
(a) the gray area is the disk $D$. 
(b) the dotted arc is $\ell'_1$ and 
the thick arc is $\ell'_2$.}
\end{figure}

\begin{Proof}
{\bf Statement (a).}
If the edge contains a crossing, 
then we can eliminate  the crossing 
on the edge by 
a C-II move. 
This contradicts the fact that the chart is minimal.

{\bf Statement (b).}
If not, 
we can eliminate the white vertex by a C-III move.

{\bf Statement (c).}
Suppose that there exists a ring $C$
such that a complementary domain of $C$ 
does not contain any white vertices.
Let $F$ be the closure of the complementary domain.
By Statement (a), 
the ring $C$ does not intersect any terminal edge 
nor free edge. 
Thus any crossing on $C$ 
is contained in a proper internal arc of $F$. 
Since $F$ is a disk or an annulus and 
since $\Gamma\cap \partial D^2=\emptyset$, 
the domain $F$ contains 
an 'innermost' disk $D$ bounded by 
an arc $\ell_1$ on $C$ and 
a proper internal arc $\ell_2$ of $F$ 
with $\partial\ell_1=\partial \ell_2$ 
such that 
for any crossing $v$ on Int $\ell_1$ 
there exists a proper internal arc 
of $D$ 
containing $v$ and a crossing on $\ell_2$
(see Fig.~\ref{fig-09}(a)). 
Let $\ell_2'$ 
be an internal arc 
with $\ell'_2\supset\ell_2$ 
such that 
$\ell_2'-\ell_2$ does not contain a crossing.
Let 
$\ell_1'$ be an arc 
outside $F$ parallel to $\ell_1$ 
with $\partial\ell_1'=\partial\ell_2'$
(see Fig.~\ref{fig-09}(b)). 
Thanks to Proposition~\ref{MoveOutFreeEdge}, 
there does not exist any black vertex 
in the disk bounded by 
$\ell_1'\cup\ell_2'$. 
Thus 
we can shift the arc $\ell_2'$ 
to the arc $\ell_1'$
by a CI-move so that 
the number of crossings decreases at least two 
(see Fig.~\ref{fig-09}(c)). 
This contradicts the fact that the chart is minimal. 
\end{Proof}

\begin{proposition}
\label{MoveOutSimpleHoop}
Let $\Gamma$ be a minimal chart in a disk $D^2$.
For any regular neighbourhood $N$ 
of $\partial D^2$ in $D^2$,
there exists a minimal chart $\Gamma'$ 
obtained from $\Gamma$ 
by C-I-M2 moves 
and ambient isotopies of $D^2$ such that
\begin{enumerate}
\item[{\rm (a)}]
$\Gamma'\cap (D^2-N)$ contains neither free edge
nor simple hoop,
\item[{\rm (b)}]
$\Gamma'\cap N$ consists of oval nests, 
simple hoops and free edges.
\end{enumerate}
\end{proposition}

\begin{Proof}
By Proposition~\ref{MoveOutFreeEdge},
we can assume that any free edge is contained 
in the regular neighbourhood $N$ 
satisfying (b).
Let $C$ be a simple hoop.
Then a complementary domain of $C$ 
does not contain any white vertices.
Let $F$ be the closure of the complementary domain.

There are two cases:
(1) $F$ is a disk,
(2) $F$ is an annulus containing $\partial D^2$.

{\bf Case 1.}
The disk $F$ does not contain any crossing.
For, if  
the disk $F$ contains a crossing, 
then 
the crossing is contained in a ring $R$ in $F$, 
because 
the hoop $C$ does not intersect any edge. 
Since $F$ does not contain any white vertices,
the ring $R$ bounds a disk without white vertices.
This contradicts Proposition~\ref{Assumption0}(c).
Hence $F$ does not contain any crossing.

Thus $\Gamma\cap F$ consists of simple hoops.
By a similar way to the proof of 
Proposition~\ref{MoveOutFreeEdge},
we can move the set $\Gamma\cap F$ 
into the neighbourhood $N$ satisfying (b) 
so that we can decrease the number of 
simple hoops in $D^2-N$. 

{\bf Case 2.}
The annulus $F$ does not contain any crossing.
For, if $F$ contains a crossing, 
then 
there exists a ring $R$ in $F$ such that 
the union $R\cup \partial D^2$ bounds 
an annulus without white vertices, 
or the ring $R$ bounds a disk in the annulus $F$ 
without white vertices. 
This contradicts 
Proposition~\ref{Assumption0}(c).
Hence $F$ does not contain any crossing.

Thus $\Gamma\cap F$ consists of oval nests, 
simple hoops and free edges.
By ambient isotopies of $D^2$,
we can move the set $\Gamma\cap F$ 
into the neighbourhood $N$ satisfying (b) 
so that we can decrease the number of 
simple hoops in $D^2-N$. 
\end{Proof}

For any minimal chart in a disk $D^2$
we can move free edges and simple hoops into 
a regular neighbourhood of $\partial D^2$ 
by C-I-M2 moves and ambient isotopies of $D^2$
by Proposition~\ref{MoveOutFreeEdge}
and Proposition~\ref{MoveOutSimpleHoop}.
Even during argument,
if free edges or simple hoops appear, 
we immediately move them 
into a regular neighbourhood of $\partial D^2$ in $D^2$.
Thus we assume the following 
(cf. \cite{OneCrossing}, \cite[assumption 1]{MinimalChart}):

\begin{assumption}
\label{AssumptionFreeEdge}
{\it For any minimal chart in a disk $D^2$, 
all the free edges and 
simple hoops  
are in a regular neighbourhood of 
$\partial D^2$ in $D^2$.}
\end{assumption}

Let $\Gamma$ be a minimal chart in 
a disk $D^2$, and 
$X$ the union of all the free edges 
and simple hoops.
Now $X$ is in a regular neighbourhood $N$ of 
$\partial D^2$ in $D^2$ 
by Assumption~\ref{AssumptionFreeEdge}.
We define\index{${\rm Main}(\Gamma)$}
$${\rm Main}(\Gamma)=\Gamma-X.$$
Let $\widehat D=Cl(D^2-N)$. 
Then $\Gamma\cap\widehat D=$Main$(\Gamma)$.
Hence $(\Gamma\cap\widehat D,\widehat D)$ is 
a tangle without free edges 
and simple hoops.

In this paper we always assume that
\begin{enumerate}
\item[] 
{\it for any tangle $(\Gamma\cap D,D)$, 
the disk $D$ does not contain 
any free edge nor a simple hoop}.
\end{enumerate}

\begin{figure}
\begin{center}
\includegraphics{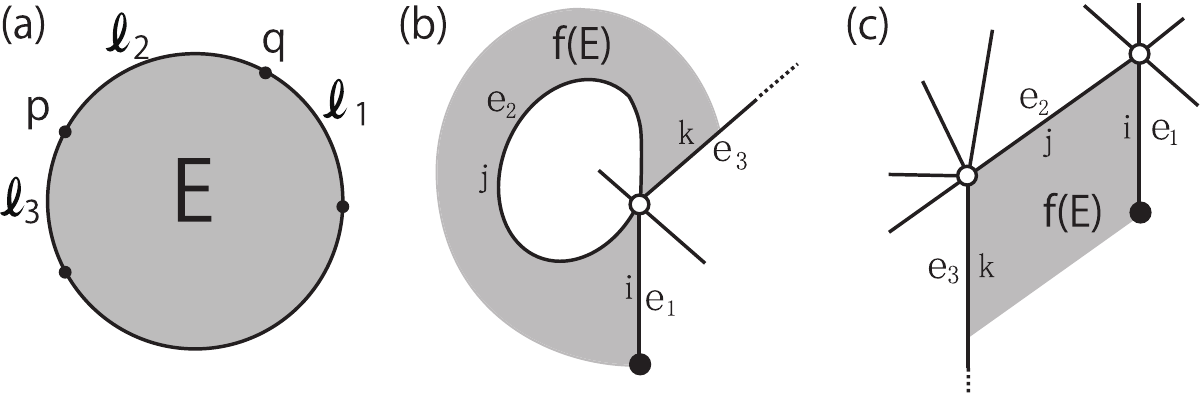}
\end{center}
\caption{ \label{fig-10} 
(b), (c) consecutive triplets.}
\end{figure}

Let $E$ be a disk, and
$\ell_1,\ell_2,\ell_3$ three arcs on $\partial E$
such that each of $\ell_1\cap \ell_2$ and $\ell_2\cap \ell_3$ is one point and $\ell_1\cap \ell_3=\emptyset$
(see Fig.~\ref{fig-10}(a)),
say $p=\ell_1\cap \ell_2$,
$q=\ell_2\cap \ell_3$.
Let $\Gamma$ be a chart in a disk $D^2$.
Let $e_1$ be a terminal edge of 
 $\Gamma$. 
A triplet $(e_1,e_2,e_3)$ of 
mutually different edges of $\Gamma$
is called 
a {\it consecutive triplet}\index{consecutive triplet}
 if there exists
a continuous map $f$ from the disk $E$ 
to the disk $D^2$ such that (see Fig.~\ref{fig-10}(b) and (c))
\begin{enumerate}
\item[(i)] the map $f$ is injective on $E-\{p,q\}$,
\item[(ii)] 
$f(\ell_3)$ is an arc in $e_3$, and $f({\rm Int}~E)\cap\Gamma=\emptyset$,
$f(\ell_1)=e_1$,
$f(\ell_2)=e_2$,
\item[(iii)]
each of $f(p)$ and $f(q)$ is a white vertex.
\end{enumerate}
If the label of $e_3$ is different
from the one of $e_1$, 
then the consecutive triplet is said to be
{\it admissible}.


\begin{remark}
{\rm Let $(e_1,e_2,e_3)$
be a consecutive triplet. 
Since $e_2$ is an edge of $\Gamma$, 
the edge $e_2$ MUST NOT contain a crossing.}
\end{remark}


\begin{lem}
\label{ConsecutiveTripletLemma} 
{\rm [Consecutive Triplet Lemma]}
{\rm $($\cite[lemma 1.1]{OneCrossing}, 
\cite[lemma 3.1]{MinimalChart}$)$}
{\it Any consecutive triplet 
in a minimal chart is admissible.}
\end{lem}

Let $\Gamma$ be a chart.
A tangle $(\Gamma\cap D,D)$ is 
called an {\it NS-tangle of label $m$} \index{NS-tangle}
(new significant tangle) 
provided that
\begin{enumerate}
\item[(i)] if $i\neq m$, 
then $\Gamma_i\cap \partial D$ is 
at most one point,
\item[(ii)] 
$\Gamma\cap D$ contains at least one white vertex, and 
\item[(iii)]
for each label $i$, 
the intersection $\Gamma_i\cap D$ contains 
at most one crossing.
\end{enumerate}

\begin{lem}\label{LemNS-Tangle}
{\rm (\cite[theorem 1.2]{MinimalChart})}
In a minimal chart, 
there does not exist 
an NS-tangle of any label.
\end{lem}

Let $\Gamma$ be a chart, and $m$ a label of $\Gamma$.
Let $E$ be a disk with $\partial E\subset \Gamma_m$.
Then the disk $E$ is called a {\it $2$-color disk}\index{2-color disk}
provided that
$\Gamma\cap E\subset 
\Gamma_{m}\cup\Gamma_{m-1}$
or 
$\Gamma\cap E\subset 
\Gamma_{m}\cup\Gamma_{m+1}$.

\begin{lem}\label{LemTwo-ColorDiskNoTerminalEdge}
{\rm (\cite[corollary 3.4(b)]{MinimalChart})}
Let $\Gamma$ be a minimal chart, 
and $m$ a label of $\Gamma$.
If $E$ is a $2$-color disk 
with $\partial E\subset\Gamma_m$, 
then $E$ does not contain any terminal edge.
\end{lem}

\begin{lem}\label{LemTwo-ColorDisk}
{\rm (\cite[lemma 7.6(b)]{MinimalChart})}
Let $\Gamma$ be a minimal chart, 
and $m$ a label of $\Gamma$.
If $E$ is a $2$-color disk 
with $\partial E\subset\Gamma_m$ 
but without free edges nor simple hoops, 
then $\Gamma_m\cap E$ is connected.
\end{lem}

Let $\Gamma$ be a chart, and 
$m$ a label of $\Gamma$. 
A simple closed curve in $\Gamma_m$ is \index{cycle}
called a {\it cycle of label $m$}.

Let $\Gamma$ be a chart,
and $m$ a label of $\Gamma$.
Let $C$ be a cycle of label $m$ in $\Gamma$ 
bounding a disk $E$.
Then an edge $e$ of label $m$ 
is called
an {\it outside edge for $C$}\index{outside edge} provided that 
\begin{enumerate}
\item[(i)]
$e\cap C$ consists of one white vertex or two white vertices, and
\item[(ii)]
$e\not\subset E$.
\end{enumerate}
For a cycle $C$ of label $m$,\index{${\mathcal{W}}_O^{{\rm Mid}}(C,m)$} 
we define
$$\begin{array}{ll}
{\mathcal{W}}(C)&= \{ w \ | \text{ $w$ is a white vertex in $C$} \},\\
{\mathcal{W}}_O^{{\rm Mid}}(C,m)&= \{ w\in \mathcal{W}(C)  \ | \text{ there exists an outside edge for $C$ {\it middle} at $w$} \}.
\end{array}$$

Let $\Gamma$ be a chart, and $m$ a label of $\Gamma$.
Let $E$ be a disk with $\partial E\subset\Gamma_m$.
Then the disk $E$ is called 
a {\it $3$-color disk}\index{3-color disk} provided that
\begin{enumerate}
\item[(i)] the disk $E$ does not contain any crossings, and
\item[(ii)] 
$\Gamma\cap E\subset \Gamma_{m-1}\cup\Gamma_m
\cup\Gamma_{m+1}$.
\end{enumerate}

\begin{lem}\label{LemThree-ColorDisk}
{\rm (\cite[corollary 4.4]{MinimalChart})}
Let $\Gamma$ be a minimal chart.
Let $C$ be a cycle of label $m$ in $\Gamma$
bounding a $3$-color disk $E$
without free edges nor simple hoops.
If $\Gamma_m\cap E$ connected,
then $|{\mathcal W}_O^{{\rm Mid}}(C,m)|\ge 2$.
\end{lem}

Combining Lemma~\ref{LemTwo-ColorDisk} 
and Lemma~\ref{LemThree-ColorDisk} 
we have the following lemma.

\begin{lem}\label{LemTwoColorTangle}
Let $\Gamma$ be a minimal chart, and 
$m,k$ labels of $\Gamma$ with $|m-k|=1$. 
Let $(\Gamma\cap D,D)$ be 
a tangle
with $\Gamma\cap D\subset\Gamma_m\cup\Gamma_k$
but without free edges nor simple hoops. 
Then for any cycle $C$ of label $m$ in $D$, 
we have 
$|{\mathcal W}_O^{{\rm Mid}}(C,m)|\ge 2$.
\end{lem}

\begin{lem}\label{BoundaryConditionLemma}
{\rm [Boundary Condition Lemma]
$($\cite[lemma 4.1]{TwoCrossingI}, 
\cite[lemma 11.1]{MinimalChart}$)$}
Let $(\Gamma\cap D,D)$ be a tangle 
in a minimal chart $\Gamma$ 
such that 
$D$ does not contain any crossing, 
free edge, nor simple hoop.
Let $a=\min\{\ i \ | \
\Gamma_i\cap\partial D\not=\emptyset\}$ and
$b=\max\{\ i \ | \
\Gamma_i\cap\partial D\not=\emptyset\}$.
Then
 $\Gamma_i\cap D=\emptyset$ 
 except for $a\le i \le b$.
\end{lem}

\section{{\large Directed Paths}}
\label{s:OneWayPath}

In this section
we investigate a path in $\Gamma_m$.

Let $\Gamma$ be a chart, 
and $m$ a label of $\Gamma$.
A simple arc $P$ in $\Gamma$
is called a {\it path}
provided that 
the end points of $P$ are 
vertices of $\Gamma$.
In particular,
if the path $P$ is in $\Gamma_m$,
then $P$ \index{path of label $m$}
is called a {\it path of label $m$}. 
Suppose that 
$v_0,v_1,\cdots,v_{p}$ 
are all the vertices in a path $P$ 
situated in this order on $P$. 
For each $i=1,\cdots,p$, 
let $e_i$ be the edge in $P$ 
with $\partial e_i=\{v_{i-1},v_{i}\}$.
Then the $(p+1)$-tuple 
$(v_0,v_1,\cdots,v_{p})$
is called 
a {\it vertex sequence}\index{vertex sequence}
of the path $P$,
and the $p$-tuple 
$(e_1,e_2,\cdots,e_{p})$
is called 
an {\it edge sequence}\index{edge sequence}
of the path $P$.

Let $P$ be a path 
in a chart with 
an edge sequence $(e_1,e_2,\dots,e_p)$.
An edge $e$ of the chart is called 
a {\it side-edge} for $P$\index{side-edge}
if $e\not\subset P$ 
but 
$e\cap
(e_2\cup e_3\cup\cdots\cup e_{p-1})
\neq\emptyset$.

Let $m$ be a label of a chart $\Gamma$, and
$P$ a path of label $m$
with 
a vertex sequence $(v_0,v_1,\cdots,v_p)$ and
an edge sequence $(e_1,e_2,\dots,e_p)$. 
The path $P$ 
is called a {\it directed path}\index{directed path}
if 
for each $i=1,2,\cdots,p$,
the edge $e_i$ is oriented 
from $v_{i-1}$ to $v_i$.
The path $P$
is called an {\it M$\&$M path}\index{M$\&$M path}
if the edge $e_1$ is middle at $v_0$ 
and 
the edge $e_p$ is middle at $v_p$.
The path $P$ is called 
a {\it dichromatic path}\index{dichromatic path} 
if there exists 
a label $k$ with $|m-k|=1$ such that
any vertex of the path 
is contained in $\Gamma_m\cap\Gamma_k$.

{\bf Warning.}
If $P$ is a directed path with 
a vertex sequence $(v_0,v_1,\cdots,v_p)$ and
an edge sequence $(e_1,e_2,\dots,e_p)$, 
then {\it we always assume that
\begin{center}
each edge $e_i~(i=1,2,\cdots,p)$ is 
oriented from $v_{i-1}$ to $v_i$.
\end{center}
}

\begin{lem}\label{LemNoMM}
In a minimal chart in a disk $D^2$, 
for any label $m$ 
there does not exist 
any dichromatic M$\&$M directed path 
of label $m$.
\end{lem}

\begin{figure}
\begin{center}
\includegraphics{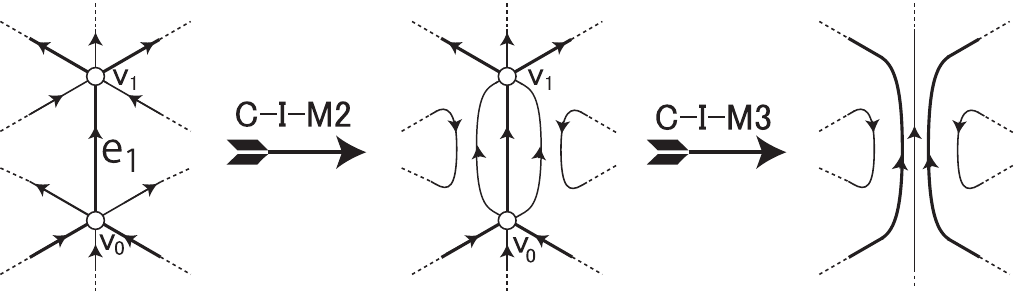}
\end{center}
\caption{ \label{fig-11}
The edge $e_1$ is middle at $v_0$ and $v_1$.}
\end{figure}

\begin{Proof}
Suppose that there exists 
a dichromatic M$\&$M directed path 
of label $m$ in
a minimal chart $\Gamma$.
Let $P$ be 
a dichromatic M$\&$M directed path of label $m$ 
containing the least number of edges 
amongst all the dichromatic 
M$\&$M directed paths 
in the chart.
Set $(v_0,v_1,\cdots,v_p)$ and 
$(e_1,e_2,\cdots,e_p)$ 
the vertex sequence and
the edge sequence of $P$.

We claim that $p=1$.
For, if $p>1$, 
let $e$ be the side-edge 
of label $m$ for $P$ 
containing the vertex $v_1$.
If the edge $e$ is outward at $v_1$, then
$e_1$ is middle at $v_1$.
The path with the edge sequence $(e_1)$ 
is a dichromatic M$\&$M directed path of label $m$  whose length is shorter than $P$.
This is a contradiction.
If the edge $e$ is inward at $v_1$, then 
$e_2$ is middle at $v_1$.
The path with the edge sequence 
$(e_2,e_3,\cdots,e_p)$ 
is a dichromatic M$\&$M directed path 
whose length is shorter than $P$.
This is a contradiction.
Thus $p=1$.

Since $e_1$ is middle at 
$v_0$ and $v_1$, 
we can 
eliminate the two vertices $v_0$ and $v_1$ 
by
two C-I-M2 moves and a C-I-M3 move 
in a regular neighbourhood of 
the edge $e_1$ in $D^2$
(see Fig.~\ref{fig-11}).
This contradicts the fact that 
the given chart is minimal.
This proves Lemma~\ref{LemNoMM}.\end{Proof}

Let $\Gamma$ be a chart, and
$m$ a label of $\Gamma$.
A simple closed curve in $\Gamma_m$ is called 
a {\it loop}\index{loop}
if it contains exactly one white vertex.
Note that loops may contain crossings of $\Gamma$.

Let $\Gamma$ be a chart, 
and $m$ a label of $\Gamma$.
Let $C$ be 
a cycle of label $m$.
Let  
$v_0,v_1,\cdots,v_{p-1}$ 
be all the vertices in $C$, 
and 
$e_1,e_2,\cdots,e_p$
all the edges in $C$.
Then the cycle $C$  
is called 
a {\it directed cycle}\index{directed cycle}
provided that 
for each $i=1,2,\cdots,p$,
the edge $e_i$ is oriented from $v_{i-1}$ to $v_i$,
where $v_p=v_0$.
We consider a loop as a directed cycle.

\begin{lem}\label{LemOneWayCycle}
Let $\Gamma$ be a minimal chart, 
$m,k$ labels of $\Gamma$ with $|m-k|=1$, and 
$(\Gamma\cap D,D)$ a tangle
with
$\Gamma\cap D\subset\Gamma_m\cup\Gamma_k$
but without free edges nor simple hoops. 
Then $D$ does not contain any directed cycle 
of label $m$.
\end{lem}

\begin{Proof}
Suppose that $D$ contains a directed cycle $C$ 
of label $m$. Now
\begin{enumerate}
\item[$(1)$]
for a white vertex $v$ in $\Gamma_m$ 
contained in 
three edges $e_1,e_2,e_3$ of label $m$, 
if $e_1$ is inward at $v$ and 
if $e_2$ is outward at $v$, 
then $e_3$ is not middle at $v$.
\end{enumerate}
Hence we have
\begin{enumerate}
\item[$(2)$] 
${\mathcal W}_O^{{\rm Mid}}(C,m)=\emptyset$.
\end{enumerate}
This contradicts 
Lemma~\ref{LemTwoColorTangle}.
Thus 
there does not exist 
a directed cycle of label $m$ in $D$.
This proves Lemma~\ref{LemOneWayCycle}.
\end{Proof}

\begin{lem}
\label{LemNoOneWayCycleIntersectNTangle}
Let $\Gamma$ be a minimal chart, and
$(\Gamma\cap D,D)$ an N-tangle or 
a net-tangle
with a label pair $(\alpha,\alpha+1)$. 
Let $m=\alpha$ or $m=\alpha+1$.
Then $D$ does not intersect 
any directed cycle of label $m$.
\end{lem}

\begin{Proof}
Suppose that $D$ intersects 
a directed cycle $C$ of label $m$.
Then
\begin{enumerate}
\item[$(1)$]
the directed cycle $C$ does not intersect $\partial D$.
\end{enumerate}
For, if $C$ intersects $\partial D$, 
then 
$C$ contains an I-edge of label $m$ 
and an O-edge of label $m$. 
On the other hand, by Condition (iii) of 
the definition of a net-tangle, 
all the edges of label $m$ 
intersecting $\partial D$ 
are I-edges or 
all the edges of label $m$ 
intersecting $\partial D$ 
are O-edges.
This is a contradiction.
Hence the directed cycle $C$ 
does not intersect $\partial D$.

Since the directed cycle $C$ intersects $D$, 
Statement $(1)$ implies that
\begin{enumerate}
\item[$(2)$]
the directed cycle $C$ is contained in $D$. 
\end{enumerate}
On the other hand,
by Condition (ii) of 
the definition of a net-tangle, we have 
$\Gamma\cap D\subset
\Gamma_\alpha\cup\Gamma_{\alpha+1}$.
Hence $D$ does not contain 
any directed cycle of label $m$
by Lemma~\ref{LemOneWayCycle}.
This contradicts $(2)$.
Therefore Lemma~\ref{LemNoOneWayCycleIntersectNTangle} holds.\end{Proof}

\begin{figure}
\begin{center}
\includegraphics{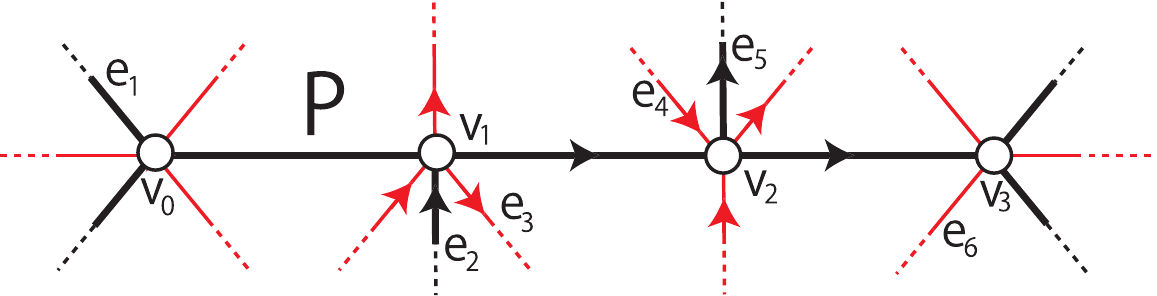}
\end{center}
\caption{ \label{fig-12}
The path $P$ is a directed path 
with a vertex sequence
$(v_0,v_1,v_2,v_3)$.}
\end{figure}

Let $P$ be a directed path of label $m$
in a chart, and  
$e$ a side-edge for $P$
not a loop.
Let $v$ be a white vertex in 
$e\cap {\rm Int~}P$, and 
$N$ a regular neighbourhood of $v$.
The edge $e$ is said to be 
{\it locally right-side }\index{locally right-side}
(resp. {\it locally left-side})\index{locally left-side} at $v$
if the arc $e\cap N$ 
is situated right (resp. left) side
of $P$ with respect to the direction of $P$ 
(see Fig.~\ref{fig-12}).
If the edge $e$ is locally 
right-side at 
a vertex $v$ 
and inward (resp. outward) at $v$, 
then
the edge 
is called a 
locally right-side edge {\it inward}\index{inward}
(resp. {\it outward})\index{outward}
at $v$. 
Similarly 
if the edge $e$ is locally 
left-side at 
a vertex $v$ 
and inward (resp. outward) at $v$, 
then
the edge 
is called a 
locally left-side edge {\it inward} 
(resp. {\it outward}) at $v$.
In Fig.~\ref{fig-12},
the edge $e_2$ is a locally right-side edge 
inward at $v_1$,
the edge $e_3$ is a locally right-side edge 
outward at $v_1$,
the edge $e_4$ is a locally left-side edge 
inward at $v_2$,
the edge $e_5$ is a locally left-side edge 
outward at $v_2$.
But $e_1$ is not 
locally left-side at $v_0$
nor $e_6$ is not 
locally right-side at $v_3$.

Let $\Gamma$ be a chart, 
and $m$ a label of the chart.
Let $e$ be an edge of label $m$, and
$P$ a directed path of label $m$ 
with an edge sequence 
$(e_1,e_2,\dots,e_p)$.
If $e=e_1$,
then the path $P$ is called 
a directed path {\it starting from} $e$.\index{starting form $e$}
Also if $e=e_p$,
then the path $P$ is called\index{leading to $e$}
a directed path {\it leading to} $e$.

\begin{figure}[tbh]
\begin{center}
\includegraphics{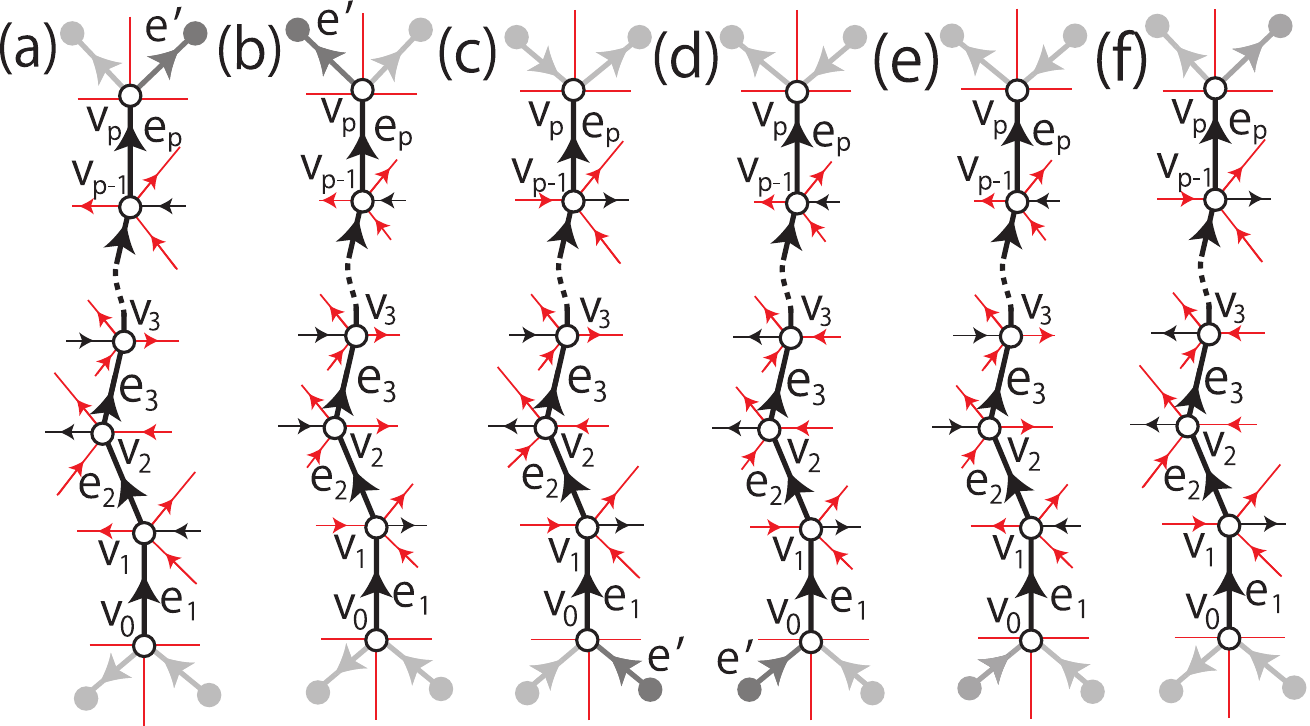}
\end{center}
\caption{\label{fig-13} 
(a) an upward-right-selective path, 
(b) an upward-left-selective path, 
(c) a downward-right-selective path, 
(d) a downward-left-selective path,
(e) an upward principal path,
(f) a downward principal path.}
\end{figure}

Let $\Gamma$ be a chart, 
and $m$ a label of the chart.
A directed path $P$ of label $m$ 
is said to be
{\it upward-right-selective} 
(resp. {\it upward-left-selective})\index{upward-right-selective}\index{upward-left-selective}
if any edge 
of label $m$
locally right-side 
(resp. left-side) 
at a vertex in Int~$P$
is inward at the vertex
(see Fig.~\ref{fig-13}(a) and (b)).
A directed path $P$ of label $m$
is said to be
{\it downward-right-selective} 
(resp. {\it downward-left-selective})\index{downward-right-selective}\index{downward-left-selective}
if any edge 
of label $m$ 
locally right-side (resp. left-side) 
at a vertex in Int~$P$
is outward at the vertex
(see Fig.~\ref{fig-13}(c) and (d)).

\begin{figure}
\begin{center}
\includegraphics{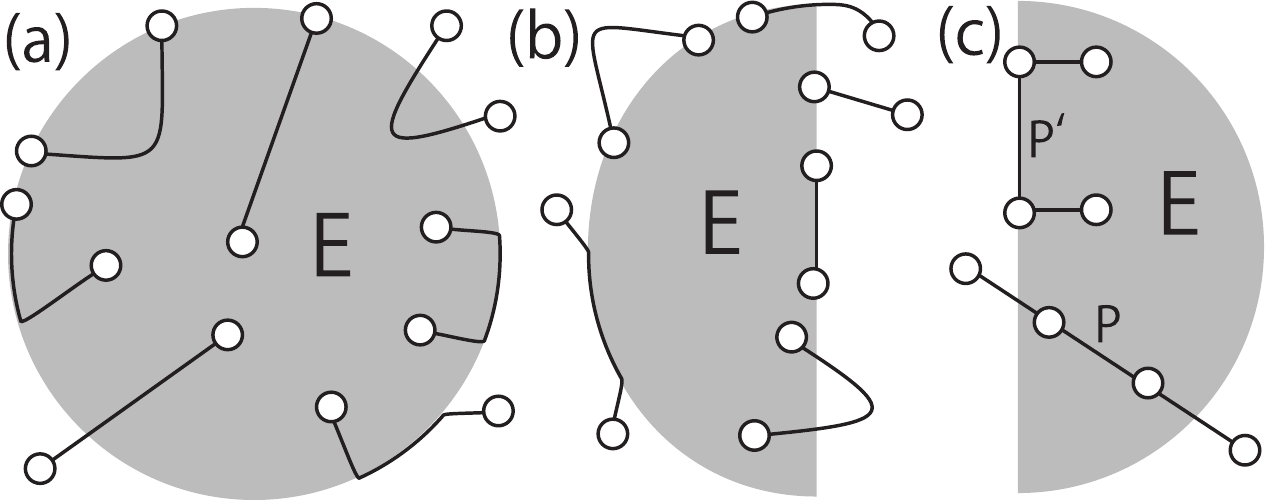}
\end{center}
\caption{\label{fig-14} 
(a) the disk $E$ dominates edges. 
(b) the disk $E$ does not dominate each edge. 
(c) the disk $E$ dominates the path $P$ but the path $P'$.}
\end{figure}

Let $\Gamma$ be a chart, 
$P$ a path, and
$E$ a disk.
If each edge in $P$ intersects 
Int~$E$ 
and if 
$P\cap E$ is connected, 
then we say that 
the path $P$ is {\it dominated by}\index{dominate} 
the disk $E$ or 
the disk $E$ {\it dominates} the path $P$ 
(see Fig.~\ref{fig-14}). 

Let $e$ be an edge 
of label $m$ in a chart, 
and $E$ a disk.
Let $P$ be an upward-right-selective 
(resp. upward-left-selective) 
directed path 
of label $m$ 
starting from $e$ 
dominated by $E$. 
The path $P$ is said to be {\it upward maximal}\index{upward maximal} 
with respect to $E$ 
if the path is not contained in another
upward-right-selective 
(resp. upward-left-selective) 
directed path of label $m$ starting from $e$
dominated by $E$. 
Let $P$ be a downward-right-selective 
(resp. downward-left-selective) 
directed path of label $m$ leading to $e$
dominated by $E$.
The path $P$ is said to be {\it downward maximal}\index{downward maximal} 
with respect to $E$ 
if the path is not contained in another 
downward-right-selective 
(resp. downward-left-selective) 
directed path of label $m$ 
leading to $e$ 
dominated by $E$.

\begin{lem}\label{lemRLselective}
Let $\Gamma$ be a minimal chart, and 
$(\Gamma\cap D,D)$ 
an N-tangle or 
a net-tangle with 
label pair $(\alpha,\alpha+1)$. 
Let $m=\alpha$ or $m=\alpha+1$. 
Let $e$ be an edge of label $m$ in $D$, 
and 
$P$ a directed path 
of label $m$ 
dominated by $D$ 
with a vertex sequence
$(v_0,v_1,\cdots,v_p)$. 
Then we have the following.
\begin{enumerate}
\item[{\rm (a)}]
Suppose that $P$ is 
an upward-right-selective 
or upward-left-selective 
directed path 
starting from $e$. 
If $v_p$ is a white vertex 
in {\rm Int}~$D$, 
then $P$ is not upward maximal with respect to $D$.
\item[{\rm (b)}]
Suppose that $P$ is 
a downward-right-selective
or downward-left-selective 
directed path 
leading to $e$. 
If $v_0$ is a white vertex 
in {\rm Int}~$D$, 
then $P$ is not downward maximal with respect to $D$. 
\end{enumerate}
\end{lem}

\begin{Proof}
{\bf Statement (a)}.
By Lemma~\ref{LemNoOneWayCycleIntersectNTangle},
there does not exist a loop of label $m$ containing $v_p$.
Now $e_p$ is inward at $v_p$.
Since $v_p$ is a white vertex, 
there are two edges of label $m$ 
at $v_p$ different from $e_p$.
There are two cases:\\
Case~$1$. The two edges are outward at $v_p$.\\
Case~$2$. One of the two edges is inward at $v_p$, 
and  the other edge is outward at $v_p$.

{\bf Case 1}. 
If $P$ is upward-right-selective 
(resp. upward-left-selective), 
let $e'$ be the one of the two edges 
outward at $v_p$
such that
the other edge is 
locally left-side 
(resp. right-side) edge at $v_p$
as a side-edge for 
the path $P\cup e'$
(we select the '{\it right}' 
(resp.'{\it left}') 
edge $e'$ with respect to 
the direction of the path $P$,
see Fig.~\ref{fig-13}(a) and (b)).

Now $P\cap e'=v_p$.
For, if $P\cap e'$ contains a vertex $v_j$
for some $0\leq j<p$, 
then 
$e_{j+1}\cup e_{j+2}\cup\cdots\cup e_p\cup e'$ 
is a directed cycle of label $m$ 
intersecting $D$.
This contradicts 
Lemma~\ref{LemNoOneWayCycleIntersectNTangle}.

Thus the path $P\cup e'$ is 
an upward-right-selective 
(resp. upward-left-selective) 
directed path of label $m$ 
starting from $e$ 
dominated by $D$. 
Hence $P$ is not upward maximal with respect to $D$.

{\bf Case 2}.
There exists only one edge 
of label $m$ 
outward at $v_p$, say $e'$.
By the similar way as the one in Case~$1$,
we can show that
if $P$ is upward-right-selective 
(resp. upward-left-selective), 
then
$P\cup e'$ is 
an upward-right-selective 
(resp. upward-left-selective) 
directed path 
starting from $e$ 
dominated by $D$. 
Hence $P$ is not upward maximal with respect to $D$.

Thus Statement (a) holds.
Similarly we can show Statement (b)
(see Fig.~\ref{fig-13}(c) and (d)).
\end{Proof}

\section{{\large N-tangles}}
\label{s:Ntangle}

In this section
we investigate a terminal edge and a directed path of label $m$ in an N-tangle.

Let $\Gamma$ be a minimal chart, 
$v$ a white vertex, 
and 
$e$ a terminal edge inward 
(resp. outward) at $v$.
Then the two edges inward (resp. outward) at $v$ 
different from $e$\index{sibling edge} 
are called the {\it sibling edges} of $e$. 

\begin{lem}\label{LemSibling}
Let $\Gamma$ be a minimal chart, 
$(\Gamma\cap D,D)$ an N-tangle or 
a net-tangle with 
a label pair $(\alpha,\alpha+1)$.
Then we have the following.
\begin{enumerate}
\item[{\rm (a)}] 
The sibling edges of 
an I-terminal edge in $D$ are O-edges.
\item[{\rm (b)}] 
The sibling edges of 
an O-terminal edge in $D$ are I-edges.
\end{enumerate}
\end{lem}

\begin{Proof}
By Condition (ii) of the definition of 
a net-tangle, we have
\begin{enumerate}
\item[$(1)$]
$\Gamma\cap D\subset\Gamma_\alpha\cup\Gamma_{\alpha+1}$.
\end{enumerate}
Let $e$ be a terminal edge in $D$,
$v$ the white vertex contained in 
the terminal edge $e$,
and $e^*$ a sibling edge of $e$. 
Suppose $e^*\subset D$. 
Since $D$ does not contain a crossing,
\begin{enumerate}
\item[$(2)$]
the edge $e^*$ does not contain a crossing.
\end{enumerate}
Since any loop is a directed cycle,
the edge $e^*$ is not a loop 
by Lemma~\ref{LemNoOneWayCycleIntersectNTangle}.
Let $v^*$ be the vertex of 
$e^*$ different from $v$.
Since $e$ is a terminal edge, 
the sibling edge $e^*$ is not a terminal edge.
Thus the vertex $v^*$ is a white vertex 
by $(2)$.
Now one of the two edges $e,e^*$ is of label $\alpha$, and
the other is of label $\alpha+1$.
Thus 
there exists a non-admissible 
consecutive triplet 
$(e,e^*,\widetilde e)$ by $(1)$ and $(2)$, 
here $\widetilde e$ 
is an edge containing the vertex $v^*$, 
and
of the same label as the one of $e$.
This contradicts
Consecutive Triplet Lemma 
(Lemma~\ref{ConsecutiveTripletLemma}).
Thus $e^*\not\subset D$.

Now the edge $e^*$ is 
an I-edge or an O-edge by
Condition (iii) of 
the definition of a net-tangle.
Therefore $e^*$ is 
an I-edge (resp. O-edge),
if $e$ is inward (resp. outward) at $v$.
Namely, 
$e^*$ is an I-edge (resp. O-edge),
if $e$ is an O-terminal (resp. I-terminal) edge.
This proves Lemma~\ref{LemSibling}.
\end{Proof}

\begin{lem}
\label{LemChromaticTangleIsN-tangle}
Let $\Gamma$ be a minimal chart, and 
$(\Gamma\cap D,D)$
a net-tangle 
with a label pair $(\alpha,\alpha+1)$. 
Then 
\begin{enumerate}
\item[{\rm (a)}]
there exists a terminal edge in $D$,
\item[{\rm (b)}]
the tangle is an N-tangle.
\end{enumerate}
\end{lem}

\begin{Proof}
By Condition (ii) of the definition of 
a net-tangle, we have
\begin{enumerate}
\item[$(1)$]
$\Gamma\cap D\subset\Gamma_\alpha\cup\Gamma_{\alpha+1}$.
\end{enumerate}
If $\Gamma\cap \partial D=\emptyset$,
then the tangle is an NS-tangle.
This contradicts Lemma~\ref{LemNS-Tangle}.
Thus $\Gamma\cap \partial D\neq\emptyset$.

Let $e_1$ be an edge 
intersecting $\partial D$.
If necessary, change orientation 
of all edges in the chart so that 
we can assume that 
\begin{enumerate}
\item[$(2)$]
the edge $e_1$ is an I-edge.
\end{enumerate}
Let $m$ be the label 
of the edge $e_1$, here 
$m=\alpha$ or $m=\alpha+1$. 
Then by Condition (iii) 
of the definition 
of a net-tangle,
\begin{enumerate}
\item[$(3)$]
all the edges of label $m$ 
intersecting $\partial D$ are I-edges.
\end{enumerate}
Let $P$ be an upward-right-selective 
directed path of label $m$ 
starting from $e_1$ upward maximal with respect to $D$.
Let $(v_0,v_1,\cdots,v_p)$ be a vertex sequence
and 
$(e_1,e_2,\cdots,e_p)$ an edge sequence of $P$.
Now the edge $e_p$ is oriented 
from $v_{p-1}$ to $v_p$.
Since there is no O-edge 
of label $m$ by $(3)$, 
we have $v_p\in$~Int~$D$.
Hence 
Lemma~\ref{lemRLselective}(a) implies 
that
the vertex $v_p$ is not a white vertex, 
because $P$ is upward maximal with respect to $D$.
Since $D$ does not contain any crossing
by Condition (i) of 
the definition of a net-tangle,
the vertex $v_p$ is a black vertex. 
Namely $e_p$ is a terminal edge 
contained in $D$.
Thus Statement (a) holds.

Since $e_p$ is inward at $v_p$, 
the edge $e_p$ is an I-terminal edge.
Thus the sibling edges of the I-terminal edge 
is O-edges by Lemma~\ref{LemSibling}(a).
Let $\beta$ be the label of the sibling edges,
here $\beta\neq m$.
Hence by Condition (iii) of the definition 
of a net-tangle,
\begin{enumerate}
\item[$(4)$]
all the edges 
of label $\beta$ intersecting 
$\partial D$ are O-edges.
\end{enumerate}
Now $(3)$ and $(4)$ imply that
if $m=\alpha$,
then the tangle is upward,
otherwise the tangle is downward.
Hence the tangle is an N-tangle.
Therefore Lemma~\ref{LemChromaticTangleIsN-tangle} holds.\end{Proof}

From now on,
for an N-tangle with 
a label pair $(\alpha,\alpha+1)$ and 
a boundary arc pair 
$(L_{\alpha},L_{\alpha+1})$,
we denote by $s_I$ 
the label of the I-edges, and 
by $s_O$ the label of the O-edges.
Further, we denote by $L_I$ (resp. $L_O$)
the one of the arcs 
$L_{\alpha},L_{\alpha+1}$
which intersects I-edges (resp. O-edges).
The pair $(s_I,s_O)$ is called\index{IO-label pair} 
the {\it IO-label pair} of the N-tangle. 
The pair $(L_I,L_O)$ is called\index{boundary IO-arc pair} 
the {\it boundary IO-arc pair} of the N-tangle. 
Then we do not need to distinguish 
upward N-tangles and downward N-tangles anymore. 
Namely,
according to our new rule,
\begin{enumerate}
\item[(I)]
any upward net-tangle 
with a label pair $(\alpha,\alpha+1)$ is
{\it an N-tangle with an IO-label pair $(s_I,s_O)$}
(here $s_I=\alpha,s_O=\alpha+1$), and
\item[(II)]
any downward net-tangle 
with a label pair $(\alpha,\alpha+1)$ is also
{\it an N-tangle with an IO-label pair $(s_I,s_O)$}
(here $s_I=\alpha+1,s_O=\alpha$).
\end{enumerate}

\begin{remark}\label{RemN-Tangle}
{\rm Let $\Gamma$ be a chart, and
$(\Gamma\cap D,D)$ an N-tangle
with an IO-label pair $(s_I,s_O)$.
\begin{enumerate}
\item[$(1)$]
There is neither O-edge of label $s_I$ nor
I-edge of label $s_O$.
\item[$(2)$]
Any terminal edge intersecting $D$ 
is contained in Int $D$
by Condition (iii) of a net-tangle and 
Condition (i) of an I-edge and an O-edge.
\item[$(3)$]
An edge is contained in Int $D$
if and only if 
its two vertices are in $D$.
\end{enumerate}
}
\end{remark}

\begin{lem}\label{LemTerminalEdgeInN-Tangle}
Let $\Gamma$ be a minimal chart, and 
$(\Gamma\cap D,D)$ an N-tangle 
with an IO-label pair $(s_I,s_O)$. 
Then we have the following.
\begin{enumerate}
\item[{\rm (a)}]
Any terminal edge of label $s_I$ in $D$ 
is an I-terminal edge.
\item[{\rm (b)}] 
Any terminal edge of label $s_O$ in $D$ 
is an O-terminal edge.
\end{enumerate}
\end{lem}

\begin{Proof} 
{\bf Statement (a)}. 
Suppose that 
there exists an O-terminal edge $e$
of label $s_I$ in $D$.
Then 
the two sibling edges of $e$ 
are of label $s_O$. 
Since any O-terminal edge is outward at 
its black vertex, 
the edge $e$ is inward at 
a white vertex $v$, and  
so are the two sibling edges.
Hence by Lemma~\ref{LemSibling}(b),
the two sibling edges are 
I-edges of label $s_O$.
This contradicts 
Remark~\ref{RemN-Tangle}$(1)$.
Hence Statement (a) holds.
Similarly we can show Statement (b).
\end{Proof}

\begin{lem}\label{LemMaximalPathFrom}
Let $\Gamma$ be a minimal chart, and 
$(\Gamma\cap D,D)$ an N-tangle
with an IO-label pair $(s_I,s_O)$. 
Let $P$ be a path dominated by $D$ 
with an edge sequence 
$(e_1,e_2,\cdots,e_{p})$. 
\begin{enumerate}
\item[{\rm (a)}] 
Suppose that $P$ is an upward-right-selective or 
upward-left-selective directed path 
starting from $e_1$. 
\begin{enumerate}
\item[{\rm (i)}] Suppose that $P$ is of label $s_I$.
Then 
$P$ is upward maximal with respect to $D$
if and only if 
$e_{p}$ is a terminal edge.
\item[{\rm (ii)}] 
Suppose that $P$ is of label $s_O$. 
Then 
$P$ is upward maximal with respect to $D$
if and only if 
$e_{p}$ is an O-edge for $D$.
\end{enumerate}
\item[{\rm (b)}]
Suppose that $P$ is 
a downward-right-selective or 
downward-left-selective directed path
leading to $e_p$.
\begin{enumerate}
\item[{\rm (i)}] 
Suppose that $P$ is of label $s_I$. 
Then 
$P$ is downward maximal with respect to $D$
if and only if 
$e_{1}$ is an I-edge for $D$.
\item[{\rm (ii)}] 
Suppose that $P$ is of label $s_O$. 
Then 
$P$ is downward maximal with respect to $D$ 
if and only if 
$e_{1}$ is a terminal edge.
\end{enumerate}
\end{enumerate}
\end{lem}

\begin{Proof}
We give a proof only 
for the case that $P$ is
upward right-selective.
Let $(v_0,v_1,\cdots,v_p)$  be 
the vertex sequence 
of the directed path $P$.

{\bf Statement (a)(i)}.
Suppose that 
the path $P$ is of label $s_I$ and
upward maximal with respect to $D$.
The edge $e_{p}$ is of label $s_I$, 
and oriented 
from $v_{p-1}$ to $v_p$.

If $v_{p}$ is outside $D$, 
then the edge $e_p$ is an O-edge of label $s_I$.
This contradicts 
Remark~\ref{RemN-Tangle}$(1)$. 
Hence the vertex 
$v_{p}$ lies in Int $D$.
Since $P$ is upward maximal with respect to $D$,
the vertex $v_p$ is not a white vertex 
by 
Lemma~\ref{lemRLselective}(a).
Since $D$ does not contain a crossing 
by Condition (i) of a net-tangle,
the edge $e_{p}$ is a terminal edge.

Conversely, if $e_p$ is a terminal edge,
then there does not exist an edge 
of label $s_I$ 
outward at $v_p$.
Hence $P$ is upward maximal with respect to $D$.

{\bf Statement (a)(ii).}
Suppose that 
the path $P$ is of label $s_O$
and upward maximal with respect to $D$.
The edge $e_p$ 
is of label $s_O$ and 
oriented 
from $v_{p-1}$ to $v_p$. 

Suppose that $v_p$ is contained in $D$.
If $v_p$ is a black vertex, then
$e_p$ is an I-terminal edge of label $s_O$. 
This contradicts 
Lemma~\ref{LemTerminalEdgeInN-Tangle}(b). 
Thus the vertex $v_p$ is 
a crossing or a white vertex.
Since $D$ does not contain a crossing 
by Condition (i) of 
the definition of a net-tangle, 
the vertex $v_p$ is a white vertex in $D$. 
Since $P$ is upward maximal with respect to $D$, 
this contradicts 
Lemma~\ref{lemRLselective}(a).
Hence $v_p$ is outside $D$. 
Since $e_p$ is oriented 
from $v_{p-1}$ to $v_p$, 
the edge $e_p$ is an O-edge.

Conversely, if $e_p$ is an O-edge,
then it is clear that $P$ is upward maximal
with respect to $D$.
\end{Proof}

\section{{\large Spindles}}
\label{s:Spindle}

In this section
we investigate  
a disk $E$ such that $\partial E$ consists of two directed paths of label $m$.

\begin{figure}
\begin{center}
\includegraphics{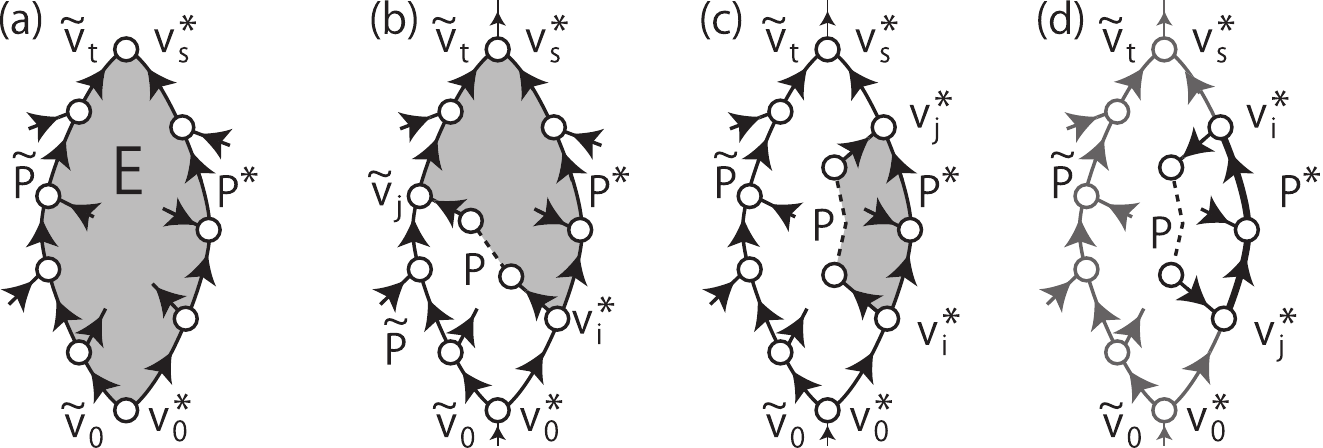}
\end{center}
\caption{\label{fig-15} Spindles.}
\end{figure}

Let $\Gamma$ be a chart, 
and $m$ a label of $\Gamma$.
Let $P^*$ be 
an upward-right-selective directed path 
of label $m$
with a vertex sequence
$(v^*_0,v^*_1,\cdots,v^*_s)$,
and $\widetilde P$ 
an upward-left-selective directed path 
of label $m$
with a vertex sequence 
$(\widetilde v_0,\widetilde v_1,\cdots,
\widetilde v_t)$.
Suppose that 
$v^*_0=\widetilde v_0$, 
$v^*_s=\widetilde v_t$, and
$P^*\cap\widetilde P=\{v^*_0,v^*_s\}$.
A disk $E$ is called a 
{\it spindle} for $\Gamma$ \index{spindle}
provided that (see Fig.~\ref{fig-15}(a))
\begin{enumerate}
\item[(i)] 
$\Gamma\cap E\subset\Gamma_m\cup\Gamma_k$
for some label $k$ with $|m-k|=1$,
\item[(ii)]
$\partial E=P^*\cup\widetilde P$, 
\item[(iii)]
for a point $x \in$Int~$P^*$, 
three points $v_0^*,x,v_s^*$
counterclockwise situate on $\partial E$
in this order.
\end{enumerate}
We call the pairs $(\widetilde P,P^*)$ and 
$(m,k)$ 
a {\it path pair} \index{path pair} and 
a {\it label pair} of the spindle $E$\index{label pair} respectively. 

\begin{remark}\label{remSpindle}
{\rm
Let $\Gamma$ be a minimal chart in a disk $D^2$, and 
$E$ a spindle for $\Gamma$
with a path pair $(\widetilde P,P^*)$ and 
a label pair $(m,k)$.
\begin{enumerate}
\item[$(1)$]
The paths $\widetilde P,P^*$ are 
dichromatic directed paths.
\item[$(2)$]
For any regular neighbourhood 
$\widehat D$ of $E$ in $D^2$, 
the pair 
$(\Gamma\cap\widehat D,\widehat D)$ 
is a tangle with
$\Gamma\cap\widehat D\subset \Gamma_m\cup\Gamma_k$. 
\item[$(3)$]
The spindle $E$ does not contain 
a directed cycle of 
label $m$ nor $k$ by $(2)$ and 
Lemma~\ref{LemOneWayCycle}.
\item[$(4)$]
Since the spindle $E$ is a $2$-color disk, 
the spindle $E$ does not contain 
a crossing nor 
a terminal edge
by the definition of a $2$-color disk and
by Lemma~\ref{LemTwo-ColorDiskNoTerminalEdge}.
\end{enumerate}
}
\end{remark}

Let $P$ be a path in a chart 
with a vertex sequence  
$(v_0,v_1,\cdots,v_p)$ 
and 
an edge sequence 
$(e_1,e_2,\cdots,e_{p})$.
For two integers $i,j$ 
with $0\le i<j\le p$, 
we denote the path 
$e_{i+1}\cup e_{i+2}\cup\cdots \cup e_{j}$\index{$P[v_i,v_j]$} 
by $P[v_i,v_j]$.

Let $\Gamma$ be a chart, 
and $E$ a disk. 
Let $e$ be an edge of $\Gamma$ 
such that 
$e\cap \partial E$ consists of 
one white vertex or two white vertices. 
If $e\subset Cl(E^c)$, then 
we call $e$  
an {\it outside edge} for $E$.\index{outside edge} 
If $e\subset E$, then 
we call $e$  
an {\it inside edge}\index{inside edge} for $E$.

\begin{lem}\label{LemProperPath}
Let $E$ be a spindle for a minimal chart $\Gamma$
with a path pair $(\widetilde P,P^*)$ and 
a label pair $(m,k)$.
Let $e$ be an inside edge for $E$
of label $m$ at
a white vertex $v$ in $\partial E$.
\begin{enumerate}
\item[{\rm (a)}]
If $e$ is outward at $v$,
then in $E$ there exists 
an upward-right-selective 
$($resp. upward-left-selective$)$ 
directed path $P$ 
of label $m$ 
starting from $e$ 
such that $P$ is a proper arc of $E$.
\item[{\rm (b)}]
If $e$ is inward at $v$,
then in $E$ there exists 
a downward-right-selective 
$($resp. downward-left-selective$)$ 
directed path $P$ 
of label $m$ 
leading to $e$ such that
 $P$ is a proper arc of $E$.
\end{enumerate}
\end{lem}

\begin{Proof}
{\bf Statement (a)}.
Let $P$ be an upward-right-selective
(resp. upward-left-selective) directed path 
of label $m$ 
starting from $e$ upward maximal 
with respect to $E$ 
with a vertex sequence
$(v_0,v_1,\cdots,v_p)$.
Suppose that $P\cap \partial E=v_0$. 
Since $E$ is a $2$-color disk, 
the disk $E$ does not contain 
a crossing nor 
a terminal edge by Remark~\ref{remSpindle}(4).
Hence $v_p$ is a white vertex 
in Int~$E$.
Since $E$ does not contain 
any directed cycle of label $m$ 
by Remark~\ref{remSpindle}(3),
we can show that $P$ is not upward maximal 
with respect to $E$
by the same way as the one 
of Lemma~\ref{lemRLselective}(a). 
This is a contradiction.
Hence $P\cap\partial E$ contains
at least two points.
Let $s=\min\{i~|~i>0,~v_i\in\partial E\}$.
Then $P[v_0,v_s]$ is a desired one.
Thus Statement (a) holds.
Similarly we can show Statement (b).
This proves Lemma~\ref{LemProperPath}.\end{Proof}

Let $\Gamma$ be a chart, and
$(\Gamma\cap D,D)$ a tangle.
For each white vertex in $D$, 
there are six short oriented arcs with
orientations inherited from 
the ones of the edges 
containing the short arcs.
Then for each white vertex, 
three consecutive arcs 
are inward at the vertex and
the other three consecutive arcs 
are outward at the vertex.
For each tangle $(\Gamma\cap D,D)$, 
we always assume that
\begin{enumerate}
\item[$(1)$]
each short arc is contained in Int~$D$, and
\item[$(2)$] 
the interiors of all the short arcs are mutually disjoint.
\end{enumerate}
Each short arc is called an {\it IS-arc} 
(resp. an {\it OS-arc})\index{IS-arc}\index{OS-arc}
{\it of the tangle}
if the short arc is inward 
(resp. outward) 
at the white vertex. 

A spindle is said to be {\it minimal}
if it does not contain another spindle.

\begin{lem} \label{lemSpindle}
For any minimal chart
there does not exist a spindle. 
\end{lem}

\begin{Proof}
Suppose that 
there exists a spindle
for a minimal chart $\Gamma$. 
Then there exists a minimal spindle $E$ 
for $\Gamma$ 
with a path pair $(\widetilde P,P^*)$ 
and a label pair $(m,k)$.
Let $(v^*_0,v^*_1,\cdots,v^*_s)$ 
and 
$(\widetilde v_0,\widetilde v_1,\cdots,
\widetilde v_t)$
be 
the vertex sequences of 
$P^*$ and $\widetilde P$ respectively.

{\bf Claim 1}.
Let $e'$ be the edge of label $m$ 
outward at $v^*_s$, and
$e''$ the edge of label $m$ 
inward at $v^*_0$. 
Then $e',e''$ 
are outside edges for $E$.

{\it Proof of Claim $1$}. 
Suppose that $e'$ is an inside edge for $E$. 
Then there exists an upward-right-selective 
directed path $P$ of label $m$ 
starting from $e'$ such that
$P$ is a proper arc of $E$. 
Hence $P\cup\widetilde P\cup P^*$ contains
a dichromatic directed cycle. 
This contradicts Remark~\ref{remSpindle}(3).
We get the same contradiction for the case 
that $e''$ is an inside edge for $E$. 
Hence Claim 1 holds.

{\bf Claim 2}. 
There does not exist an inside edge for $E$ 
of label $m$ outward at
a vertex in Int~$P^*\cup$~Int~$\widetilde P$.

{\it Proof of Claim $2$}.
Suppose that there exists an inside edge 
$\widehat e$ for $E$ of label $m$ 
outward at $v^*_i$
for some $i\in\{1,2,\cdots,s-1\}$.
By Lemma~\ref{LemProperPath}(a) 
there exists 
an upward-left-selective directed path
$P$ of label $m$ 
starting from $\widehat e$ 
with $P\cap\partial E$ 
two white vertices.
Let $v$ be the one of the two white vertices
different from $v^*_i$.
Then $v\in\{v^*_1,v^*_2,\cdots,v^*_{s-1}\}
\cup\{\widetilde v_1,\widetilde v_2,
\cdots,\widetilde v_{t-1}\}$ 
by Claim~$1$.

Suppose that $v=\widetilde v_j$ for some 
$j\in \{1,2,\cdots,t-1\}$.
Then we can find a spindle in $E$
with a path pair  
$( P \cup 
\widetilde P[\widetilde v_j,\widetilde v_t],
P^*[v^*_i,v^*_s])$
(see Fig.~\ref{fig-15}(b)). 
This contradicts the fact that 
$E$ is a minimal spindle.

Suppose that 
$v=v^*_j$ for some 
$j\in \{1,2,\cdots, {s-1}\}$.
Then $i<j$ or $j<i$.
If $i<j$, then
we can find a spindle in $E$ 
with a path pair $(P,P^*[v^*_i,v^*_j])$
(see Fig.~\ref{fig-15}(c)).
This contradicts the fact that 
$E$ is a minimal spindle.
If $j<i$, 
then $P^*[v^*_j,v^*_i]\cup P$ is 
a dichromatic directed cycle
(see Fig.~\ref{fig-15}(d)).
This contradicts Remark~\ref{remSpindle}(3).
Hence the edge $\widehat e$ is not 
outward at a vertex in Int~$P^*$.
Similarly we can show that 
there does not exist an inside edge for $E$
of label $m$ 
outward at 
a vertex in Int~$\widetilde P$.
Therefore Claim~$2$ holds.

{\bf Claim 3}.
There does not exist an inside edge for $E$
of label $m$ inward at 
a vertex in Int~$P^*\cup$~Int~$\widetilde P$.

{\it Proof of Claim~$3$}.
Suppose that there exists an inside edge 
$\widehat e$ for $E$
of label $m$ inward at $\widehat v$
in Int~$P^*\cup$~Int~$\widetilde P$.
By Lemma~\ref{LemProperPath}(b) 
there exists a directed path
$P$ of label $m$ 
leading to $\widehat e$ 
with 
$P\cap\partial E$ 
two white vertices. 
Let $v$ be 
the one of the two white vertices
different from $\widehat v$.
Thus there exists an inside edge for $E$ 
of label $m$ 
outward at $v\in$Int~$P^*\cup$~Int~$\widetilde P$.
This contradicts Claim~$2$.
Hence 
the edge $\widehat e$ 
is not inward at 
a vertex in Int~$P^*\cup$~Int~$\widetilde P$.
Hence Claim~$3$ holds.

{\bf Claim 4}.
There does not exist any white vertex in 
Int~$E$.

For, if there exists a white vertex 
in Int~$E$,
let $A$ be a regular neighbourhood of 
$\partial E$ in $E$.
Then $(A-\partial E)\cap \Gamma_m=\emptyset$ 
by Claim~$1$, Claim~$2$ and Claim~$3$.
Set $D'=Cl(E-A)$. 
Then the pair $(\Gamma\cap D',D')$ is 
an NS-tangle of label $k$.
This contradicts Lemma~\ref{LemNS-Tangle}.
Hence Claim~$4$ holds.

Now,
Claim~$2$ and Claim~$3$ imply that
\begin{enumerate}
\item[$(1)$]
any side-edge of label $m$ for 
$P^*$ or $\widetilde P$ lies 
outside $E$.
\end{enumerate}
Since there is no terminal edge nor crossing 
in $E$
by Remark~\ref{remSpindle}$(4)$,
each edge of label $k$ in $E$ contains
exactly one IS-arc and one OS-arc.
Hence
\begin{enumerate}
\item[$(2)$] 
the number of IS-arcs of label $k$ in $E$
equals 
the number of OS-arcs of label $k$ in $E$.
\end{enumerate}
Since $P^*$ is upward-right-selective,
and
since $\widetilde P$ is upward-left-selective,
Statement $(1)$ implies that 
any side-edge of label $m$ for 
$P^*$ or $\widetilde P$ is
inward at a vertex in 
$\partial E=P^*\cup\widetilde P$.
Namely for each vertex in $\partial E-v^*_s$,
the disk $E$ contains 
exactly one edge of label $k$
outward at the vertex.
Hence considering $v^*_0=\widetilde v_0$, 
the four claims Claim~$1$,  
Claim~$2$, Claim~$3$, Claim~$4$ imply that
the number of OS-arcs of label $k$ in $E$
is $s+t-1$.
On the other hand, 
there exists 
only one IS-arc of label $k$ in $E$,
which is inward at $v^*_s=\widetilde v_t$.
Thus Statement $(2)$ implies that
$s+t-1=1$. 
Since $P^*$ and $\widetilde P$ are paths, 
we have $s\ge 1$ and $t\ge 1$.
Hence $s=t=1$.
Thus $P^*$ (resp. $\widetilde P$) consists of 
exactly one edge $e^*$ (resp. $\widetilde e$). 
Further, $E$ contains an edge $e$ of label $k$
connecting the two white vertices 
$v^*_0$ and $v^*_1$ so that 
$e^*,e,\widetilde e$ are 
three consecutive edges connecting 
the two vertices $v_0^*,v_1^*$.
Hence we can eliminate the two white vertices
by a C-I-M3 move. 
This contradicts the fact that
the chart is minimal.
Therefore Lemma~\ref{lemSpindle} holds.
\end{Proof}

\section{{\large Half Spindles}}
\label{s:HalfSpindle}

In this section,
we investigate  
a disk $E$ such that $\partial E$ consists of an arc and two directed paths of label $m$.

\begin{figure}
\begin{center}
\includegraphics{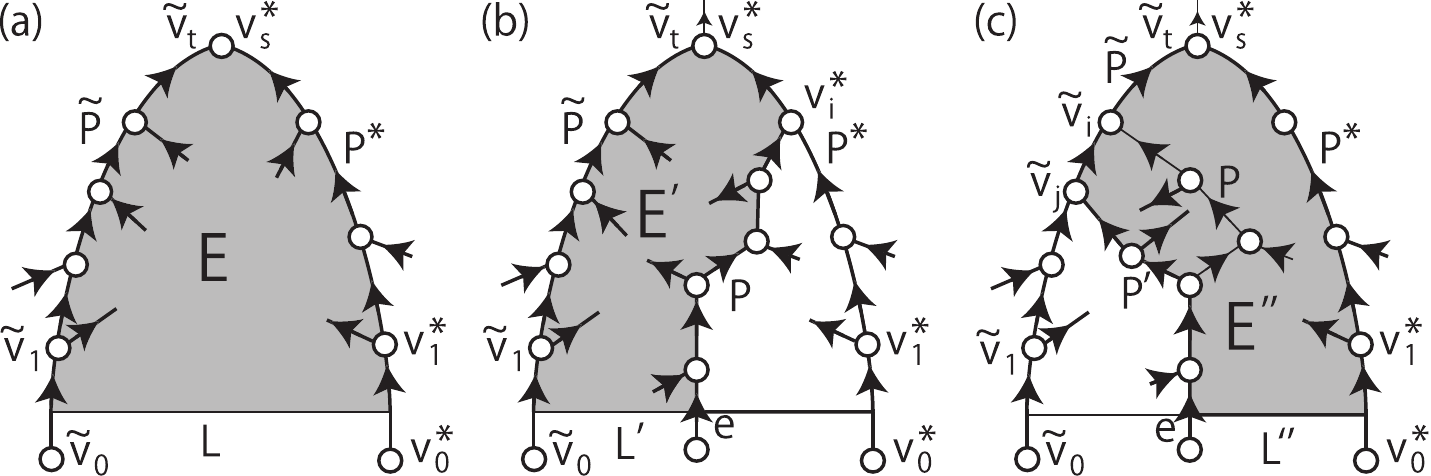}
\end{center}
\caption{\label{fig-16} Half spindles.}
\end{figure}

Let $\Gamma$ be a chart, 
and $m$ a label of $\Gamma$.
Let $P^*$ be 
an upward-right-selective directed path 
of label $m$
with a vertex sequence
$(v^*_0,v^*_1,\cdots,v^*_s)$,
and $\widetilde P$ 
an upward-left-selective directed path 
of label $m$
with a vertex sequence 
$(\widetilde v_0,\widetilde v_1,\cdots,
\widetilde v_t)$
with 
$P^*\cap \widetilde P=v^*_s=\widetilde v_t$.
Let $E$ be a disk 
and $L$ an arc on $\partial E$.
The disk $E$ is called a 
{\it half spindle} for $\Gamma$\index{half spindle}
provided that 
(see Fig.~\ref{fig-16}(a))
\begin{enumerate}
\item[(i)]
$\Gamma\cap E\subset\Gamma_m\cup\Gamma_k$
for some label $k$ with $|m-k|=1$,
\item[(ii)] 
$(P^*\cup\widetilde P)\cap$~Int~$E=\emptyset$,
and
$E\not\ni v^*_0,\widetilde v_0$, 
\item[(iii)]
$P^*\cap E$ 
is an arc containing 
$P^*[v^*_1,v^*_s]$,
$\widetilde P\cap E$
is an arc containing 
$\widetilde P[\widetilde v_1,\widetilde v_t]$,
and
$L=
Cl(\partial E-(P^*\cup\widetilde P))$,
\item[(iv)] 
the three arcs 
$\widetilde P\cap E,L,P^*\cap E$ 
are counterclockwise situated 
on $\partial E$ in this order,
\item[(v)]
if an edge 
intersects Int~$L$,
then it is an I-edge of label $m$ for $E$, 
\item[(vi)]
$s\ge 2$ and $t\ge 2$.
\end{enumerate}
We call
the triplet 
$(\widetilde P,L,P^*)$ and
the pair $(m,k)$ 
an {\it arc triplet} and
a {\it label pair} of the half spindle $E$
respectively.

\begin{remark}\label{remHalfSpindle}
{\rm
Let $\Gamma$ be a minimal chart in a disk $D^2$,
and $E$ a half spindle 
for  $\Gamma$
with an arc triplet 
$(\widetilde P,L,P^*)$ and 
a label pair $(m,k)$.
\begin{enumerate}
\item[$(1)$]
The paths 
$\widetilde P[\widetilde v_1,\widetilde v_t],
P^*[v^*_1,v^*_s]$ are 
dichromatic directed paths.
\item[$(2)$]
For each regular neighbourhood $\widehat D$ 
of $E$ in $D^2$,
the pair 
$(\Gamma\cap\widehat D,\widehat D)$ 
is a tangle with
$\Gamma\cap\widehat D
\subset \Gamma_m\cup\Gamma_k$. 
\item[$(3)$]
The half spindle $E$ 
does not contain a directed cycle of 
label $m$ nor $k$ by $(2)$ and 
Lemma~\ref{LemOneWayCycle}.
\item[$(4)$]
The half spindle $E$ does not contain 
a crossing nor 
a terminal edge 
of label $m$.
For, there does not exist
an I-edge nor an O-edge for $E$ of label $k$ by
Condition (v) of the definition of 
a half spindle.
If $E$ contains a terminal edge 
of label $m$, 
we can find a non-admissible consecutive triplet 
that contradicts 
Lemma~\ref{ConsecutiveTripletLemma}.
\end{enumerate}
}
\end{remark}

We can show the following lemma by 
a similar argument as the one of 
Lemma~\ref{LemProperPath}.

\begin{lem}\label{LemProperPathInHalfSpindle}
Let $E$ be a half spindle 
for a minimal chart $\Gamma$
with an arc triplet 
$(\widetilde P,L,P^*)$ and 
a label pair $(m,k)$.
Let $e$ be an edge of label $m$ with 
$e\cap{\rm Int}~E\not=\emptyset$ and
$e\cap\partial E$ a point.
\begin{enumerate}
\item[{\rm (a)}]
If $e$ is outward at a white vertex on $\partial E$
or if $e$ intersects $L$,
then there exists 
an upward-right-selective 
$($resp. upward-left-selective$)$ 
directed path $P$ 
of label $m$ 
starting from $e$ 
dominated by $E$ 
with 
$P\cap \partial E$ two points.
\item[{\rm (b)}]
If $e$ is inward at a white vertex on $\partial E$,
then there exists a downward-right-selective
 $($resp. downward-left-selective$)$ 
directed path $P$ 
of label $m$ 
leading to $e$ 
dominated by $E$ 
with 
$P\cap \partial E$ two points.
\end{enumerate}
\end{lem}

A half spindle is said to be {\it minimal}
if it does not contain another half spindle.

\begin{lem}\label{lemHalfSpindle}
For any minimal chart 
there does not exist a half spindle. 
\end{lem}

\begin{Proof}
Suppose that there exists a half spindle
for a minimal chart $\Gamma$. 
Then there exists a minimal half spindle $E$
with an arc triplet 
$(\widetilde P,L,P^*)$
and a label pair $(m,k)$. 
Let $(v^*_0,v^*_1,\cdots,v^*_s)$ 
and 
$(\widetilde v_0,\widetilde v_1,\cdots,
\widetilde v_t)$
be the vertex sequences of 
$P^*$ and $\widetilde P$
respectively.

{\bf Claim 1}.
Let $e'$ be the edge of label $m$ outward at $v^*_s$.
Then $e'$ is an outside edge for $E$. 

{\it Proof of Clain~$1$}.
Suppose that $e'$ is an inside edge for $E$.
Then by Lemma~\ref{LemProperPathInHalfSpindle}(a),
there exists an upward-right-selective
directed path $P$ of label $m$ starting from $e$
dominated by $E$ 
with $P\cap\partial E$ two points.
Since there does not exist an O-edge 
intersecting ${\rm Int}~L$ by Condition (v)
of the definition of a half spindle, 
we have $P\cap{\rm Int}~L=\emptyset$.
Thus
$P\cap\partial E=P\cap (\widetilde P \cup P^*)$
consists of two white vertices. 
Hence $P\cup \widetilde P \cup P^*$ 
contains a directed cycle.
This contradicts Remark~\ref{remHalfSpindle}(3).
Hence $e'$ is an outside edge for $E$. 
Therefore Claim~$1$ holds.

{\bf Claim 2}.
$\Gamma_m\cap{\rm Int}~L=\emptyset.$

{\it Proof of Claim~$2$}.
Suppose 
$\Gamma_m\cap{\rm Int}~L
\neq\emptyset.$ 
Let $e$ be an edge of label $m$
intersecting Int~$L$.
By Lemma~\ref{LemProperPathInHalfSpindle}(a), 
there exists an upward-right-selective 
directed path $P$ of label $m$ 
starting from $e$ with
$P\cap\partial E$ two points.
Since there does not exist 
an O-edge intersecting Int~$L$, 
the intersection
$P\cap(P^*\cup\widetilde P)$ is a white vertex $v$.

If $v=v^*_i~(1\leq i<s)$, 
then $P$ splits the disk $E$ 
into two disks. 
Let $E'$ be the one of the two disks
intersecting $\widetilde P$.
Then $E'$ is a half spindle with 
an arc triplet
$(\widetilde P,L',P\cup P^*[v^*_i,v^*_s])$
for some arc $L'\subset L$
(see Fig.~\ref{fig-16}(b)).
This contradicts the fact that 
$E$ is a minimal half spindle.

Suppose $v=\widetilde v_i~(1\le i<t)$.
Then $P$ splits the disk $E$ into 
two disks.
Let $E'$ be the one of the two disks
not intersecting $P^*$.
Let $P'$ be 
an upward-left-selective directed path
of label $m$ starting from $e$
with 
$P'\cap(P^*\cup\widetilde P)$ 
a white vertex $v'$.
Then $P'\subset E'$ and 
hence $v'=\widetilde v_j~(1\le j\le i<t)$.
Thus $P'$ splits $E$ into two disks.
Let $E''$ be 
the one of the two disks 
intersecting $P^*$.
Then $E''$ is a half spindle with an arc triplet
$(P'\cup\widetilde P[\widetilde v_j,\widetilde v_t],
L'',P^*)$ for some arc $L''\subset L$
(see Fig.~\ref{fig-16}(c)).
This contradicts the fact 
that $E$ is a minimal half spindle.
Therefore Claim~$2$ holds.

By the help of Claim~$2$, 
we can show the following  
Claim~$3$ and Claim~$4$ 
by the same argument 
as the ones of
Claim~$2$ and Claim~$3$ of 
Lemma~\ref{lemSpindle} 
respectively.

{\bf Claim 3}.
There does not exist an inside edge for $E$ 
of label $m$ outward at
a vertex in 
Int~$P^*\cup{\rm Int}~\widetilde P$.

{\bf Claim 4}.
There does not exist an inside edge for $E$ 
of label $m$ inward at
a vertex in 
Int~$P^*\cup{\rm Int}~\widetilde P$.

By a similar argument as the one of Claim~$4$ 
in Lemma~\ref{lemSpindle},
the four claims Claim~$1$, Claim~$2$,
Claim~$3$ and Claim~$4$ assure 
that
\begin{enumerate}
\item[$(1)$] there does not exist 
any white vertex in 
Int~$E$.
\end{enumerate}
Since $s\ge2,t\ge2$ by Condition (vi) of the definition of a half spindle,
Lemma~\ref{ConsecutiveTripletLemma} and
 Claim~$2$ imply that 
there does not exist 
a terminal edge 
of label $k$ in $E$.
Thus 
each edge of label $k$ in $E$ 
contains
exactly one IS-arc and one OS-arc.
Hence
\begin{enumerate}
\item[$(2)$] 
the number of IS-arcs of label $k$ in $E$
equals 
the one of OS-arcs of label $k$ in $E$.
\end{enumerate}
Now,
Claim~$3$ and Claim~$4$ imply that
\begin{enumerate}
\item[$(3)$] any side-edge of label $m$ for 
$P^*$ or $\widetilde P$ lies 
outside $E$.
\end{enumerate}
Since $P^*$ is 
upward-right-selective,
and
since $\widetilde P$ is 
upward-left-selective,
Statement $(3)$ implies that 
any side-edge of label $m$ for 
$P^*$ or $\widetilde P$ is
inward at a white vertex in 
$P^*\cup\widetilde P$.
Thus for each white vertex 
in $\partial E-v^*_s$,
the disk $E$ contains 
exactly one edge of label $k$
outward at the vertex.
Hence considering 
$E\not\ni v^*_0,\widetilde v_0$, 
the three claims Claim~$1$, 
Claim~$3$, Claim~$4$, 
and Statement $(1)$ imply that
the number of OS-arcs of label $k$ in $E$
is $s+t-2$.
On the other hand, 
there exists 
only one IS-arc of label $k$ in $E$,
which is inward at 
$v^*_s=\widetilde v_t$.
Thus Statement $(2)$ implies that
$s+t-2=1$, i.e. $s+t=3$. 
On the other hand, 
by Condition (vi) 
of the definition 
of a half spindle, 
we have $s+t\ge 2+2=4$.
This is a contradiction.
Hence Lemma~\ref{lemHalfSpindle} holds.
\end{Proof}

\section{{\large Proof of Theorem~1.1}}
\label{s:ProofTh1}

In this section
we prove Theorem~\ref{Theorem1}.

Let $\Gamma$ be a chart.
For each tangle $(\Gamma\cap D,D)$, 
we define

$E_I(D)=$ the number of I-edges for $D$,\index{$E_I(D)$}

$E_O(D)=$ the number of O-edges for $D$,\index{$E_O(D)$}

$T_I(D)=$ the number of I-terminal edges\index{$T_I(D)$}
in $D$,\index{$T_O(D)$}

$T_O(D)=$ the number of O-terminal edges 
in $D$.

\begin{lem}\label{EI+TI=EO+TO}
Let $\Gamma$ be 
a minimal chart, 
and 
$(\Gamma\cap D,D)$ a tangle without crossing. 
If $\partial D$ does not intersect 
any terminal edge,
then we have
\center{$E_I(D)+T_O(D)=E_O(D)+T_I(D)$}.
\end{lem}

\begin{Proof}
Since $\partial D$ does not intersect any terminal edge,
any terminal edge intersecting $D$
is contained in $D$.
There are three IS-arcs and three OS-arcs
around each white vertex in $D$. 
Hence
\begin{enumerate}
\item[$(1)$] 
the number of IS-arcs in $D$ is equal to 
the number of OS-arcs in $D$.
\end{enumerate}
Let 
$\Bbb S=
\{e~|~e$ is an edge of $\Gamma$
containing two white vertices in $D\}$, and\\
$\Bbb L=
\{\ell~|~\ell$ is a loop in $\Gamma$ 
containing a vertex in $D\}$. 
Then we have the following.
\begin{enumerate}
\item[$(2)$]
Each I-edge contains exactly one IS-arc
in $D$.
\item[$(3)$]
Each O-edge contains exactly one OS-arc
in $D$.
\item[$(4)$] 
Each I-terminal edge contains 
exactly one OS-arc.
\item[$(5)$] 
Each O-terminal edge contains
exactly one IS-arc.
\item[$(6)$] 
Each edge in $\Bbb S$ contains 
exactly one IS-arc and 
exactly one OS-arc. 
\item[$(7)$] 
Each loop in $\Bbb L$ contains 
exactly one IS-arc and 
exactly one OS-arc. 
\end{enumerate}
Thus Statement $(1)$ implies 
\begin{enumerate}
\item[] 
$E_I(D)+T_O(D)+|\Bbb S|+|\Bbb L|
=E_O(D)+T_I(D)+|\Bbb S|+|\Bbb L|$.
\end{enumerate}
Therefore 
$E_I(D)+T_O(D)=E_O(D)+T_I(D)$. 
This proves Lemma~\ref{EI+TI=EO+TO}.\end{Proof}

Let $P$ be a directed path of label $m$ 
in a chart with 
a vertex sequence $(v_0,v_1,\cdots,v_p)$ and 
an edge sequence $(e_1,e_2,\cdots,e_p)$. 
The path $P$ is {\it upward principal} 
provided that for each  $i=1,2,\cdots,p$ 
the edge $e_i$ is middle at $v_{i-1}$. 
The path $P$ is {\it downward principal} 
provided that for each  $i=1,2,\cdots,p$ 
the edge $e_i$ is middle at $v_{i}$. 

\begin{remark}
\label{remarkPrincipal}
{\rm 
Let $\Gamma$ be a chart, and 
$m$ a label of $\Gamma$.
Let $P$ be a directed path 
of label $m$ in $\Gamma$ 
with a vertex sequence 
$(v_0,v_1,\cdots,v_p)$ and 
an edge sequence 
$(e_1,e_2,\cdots,e_p)$. 
\begin{enumerate}
\item[(1)]
If $P$ is upward principal,
then  
for any edge $e$ of label $m$ with 
$e\cap P=v_i$ for some $i~(0<i<p)$, 
the edge $e$ is inward at $v_i$ 
(see Fig.~\ref{fig-13}(e)).
\item[(2)] 
If $P$ is upward principal, 
then 
$P$ is upward-right-selective and  
upward-left-selective.
\item[(3)] 
If $P$ is downward principal 
for any edge $e$ of label $m$ with 
$e\cap P=v_i$ for some $i~(0<i<p)$, 
then 
the edge $e$ is outward at $v_i$ 
(see Fig.~\ref{fig-13}(f)).
\item[(4)] 
If $P$ is downward principal, 
then 
$P$ is downward-right-selective and  
downward-left-selective.
\end{enumerate}
} 
\end{remark}

\begin{lem}
\label{lemPrincipalPath}
Let $\Gamma$ be a minimal chart, and 
$m$ a label of $\Gamma$.
Let $P$ be a dichromatic directed path 
of label $m$ in $\Gamma$ 
with a vertex sequence 
$(v_0,v_1,\cdots,v_p)$ and 
an edge sequence 
$(e_1,e_2,\cdots,e_p)$. 
\begin{enumerate}
\item[{\rm (a)}] 
If $e_1$ is middle at $v_0$, 
then $P$ is upward principal.
\item[{\rm (b)}] 
If $e_p$ is middle at $v_p$, 
then $P$ is downward principal.
\end{enumerate}
\end{lem}

\begin{Proof}
Suppose that $e_1$ is middle at $v_0$. 
Let $i$ be an integer with $0<i<p$, and 
$e$ an edge of label $m$ with 
$e\cap P=v_i$. 
If $e$ is outward at $v_i$, then 
$e_i$ is middle at $v_i$ 
because $e_{i+1}$ is outward at $v_i$ 
by the warning for a directed path. 
Then $P[v_0,v_i]$ is a dichromatic M\&M 
directed path. 
This contradicts Lemma~\ref{LemNoMM}. 
Thus $e$ is inward at $v_i$. 
Hence $e_{i+1}$ is middle at $v_i$. 
Thus $P$ is upward principal.
Hence Statement (a) holds. 
Similarly we can show Statement (b). 
\end{Proof}

Let $\Gamma$ be a minimal chart, and 
$(\Gamma\cap D,D)$ an N-tangle with
an IO-label pair $(s_I,s_O)$.
Let $P$ be a directed path of label $s_I$ with 
a vertex sequence $(v_0,v_1,\cdots,v_p)$
and
an edge sequence $(e_1,e_2,\cdots,e_p)$
such that
\begin{enumerate}
\item[(i)] $P\subset D$,
\item[(ii)] the vertex $v_0$ is 
a white vertex of 
a terminal edge $\tau_O$ of label $s_O$,
and
\item[(iii)] the edge $e_p$ is 
a terminal edge $\tau_I$.
\end{enumerate}
Then $e_1$ is middle at $v_0$.
Hence $P$ is upward principal by
Lemma~\ref{lemPrincipalPath}(a).
The path $P$ is
called the 
{\it upward principal path 
associated to 
a terminal edge $\tau_O$}.\index{principal path}
The terminal edge $\tau_I$ is called 
a {\it corresponding terminal edge for $\tau_O$}. 
We also say that 
the path $P$ is 
the upward principal path 
{\it connecting 
the terminal edge $\tau_O$ 
and 
the terminal edge $\tau_I$}.

\begin{lem}\label{LemPrincipalAssociated}
Let $\Gamma$ be a minimal chart, and 
$(\Gamma\cap D,D)$ an N-tangle with
an IO-label pair $(s_I,s_O)$.
Let $\tau_O$ be 
a terminal edge of label $s_O$ in $D$.
Then there exists an upward principal path 
of label $s_I$ 
connecting the terminal edge $\tau_O$ and
a terminal edge of label $s_I$.
\end{lem}

\begin{Proof}
Let $v_0$ be 
the white vertex of $\tau_O$, and 
$e_1$ the edge of label $s_I$ 
middle at $v_0$.
Since $\tau_O$ is inward at $v_0$ by 
Lemma~\ref{LemTerminalEdgeInN-Tangle}(b), 
the edge $e_1$ is outward at $v_0$.
Let $P$ be 
an upward-right-selective directed path 
of label $s_I$
starting from $e_1$ 
upward maximal with respect to $D$.
Let 
$(e_1,e_2,\cdots,e_p)$ be 
an edge sequence of $P$.
Then $e_p$ is a terminal edge 
by Lemma~\ref{LemMaximalPathFrom}(a)(i).
Therefore $P$ is an upward principal path 
of label $s_I$
associated to the terminal edge $\tau_O$, and
$e_p$ is a corresponding terminal edge for $\tau_O$.
This proves Lemma~\ref{LemPrincipalAssociated}.
\end{Proof}

\begin{remark}
\label{RemPricipalAssociated}
{\rm 
In the proof of the lemma above,
the edge $e_1$ is uniquely determined by 
the terminal edge $\tau_O$.
Hence in 
Lemma~\ref{LemPrincipalAssociated} above,
$\tau_O$ uniquely determines 
the upward-right-selective 
directed path starting from $e_1$ 
upward maximal with respect to $D$
and  
the corresponding terminal edge 
for $\tau_O$.
}
\end{remark}

Let $\Gamma$ be a minimal chart, and 
$(\Gamma\cap D,D)$ an N-tangle with
an IO-label pair $(s_I,s_O)$. 
Let 

$\mathbb T_O(D)=$ 
the set of all the terminal edges 
of label $s_O$ in $D$, 
and\index{$\mathbb T_O(D)$}

$\mathbb T_I(D)=$\index{$\mathbb T_I(D)$} 
the set of all the terminal edges 
of label $s_I$ in $D$.

By Lemma~\ref{LemPrincipalAssociated} and 
Remark~\ref{RemPricipalAssociated},
we define a map $f_D:\mathbb T_O(D)\to\mathbb T_I(D)$ 
by, for each $\tau_O\in\mathbb T_O(D)$,  
\begin{enumerate}
\item[]
$f_D(\tau_O)=$ 
the corresponding terminal edge 
for $\tau_O$.
\end{enumerate}

Let $\Gamma$ be a minimal chart, and 
$(\Gamma\cap D,D)$ an N-tangle
with an IO-label pair $(s_I,s_O)$.
For a terminal edge $\tau$ in $D$ 
with sibling edges $e^*,e^{**}$, 
the union $e^*\cup e^{**}$ splits the disk $D$
into two disks 
by Lemma~\ref{LemSibling}. 
We denote by $\Delta(\tau)$ 
the one of the two disks 
containing the terminal edge $\tau$.

\begin{lem}
\label{lemDelta}
Let $\Gamma$ be a minimal chart, and 
$(\Gamma\cap D,D)$ an N-tangle 
with an IO-label pair $(s_I,s_O)$ 
and a boundary IO-arc pair 
$(L_I,L_O)$.
Let $\tau$ be a terminal edge in $D$. 
Then we have the following.
\begin{enumerate}
\item[{\rm (a)}] 
If the label of $\tau$ is $s_O$, 
then $\Delta(\tau)\cap\partial D\subset L_I$
otherwise 
$\Delta(\tau)\cap\partial D\subset L_O$.
\item[{\rm (b)}] 
$\Gamma\cap ({\rm Int}~\Delta(\tau)-\tau)
=\emptyset$.
\end{enumerate}
\end{lem}

\begin{figure}
\begin{center}
\includegraphics{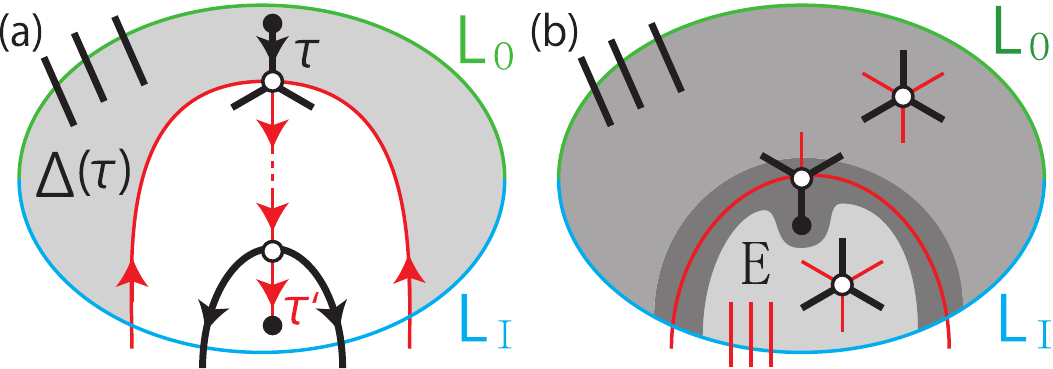}
\end{center}
\caption{\label{fig-17} 
(a) the light gray disk is $\Delta(\tau)$. (b) the light gray disk is $E$.}
\end{figure}

\begin{Proof} 
Let $e^*,e^{**}$ be the sibling edges of $\tau$, 
and
$L=\Delta(\tau)\cap\partial D$. 

{\bf Statement (a)}.
Suppose that the label of $\tau$ is $s_O$.
Since $e^*,e^{**}$ are of label $s_I$, 
the sibling edges $e^*,e^{**}$ of $\tau$ intersect 
$L_I$. 
Hence 
$L\subset L_I$ or $L\supset L_O$.

If $L\supset L_O$, 
then $(D-\Delta(\tau))\cap\partial D\subset L_I$.
Now $Cl(D-\Delta(\tau))$ contains 
the upward principal path associated to $\tau$ and 
the corresponding terminal edge $\tau'$ for $\tau$
(see Fig.~\ref{fig-17}(a)).
Since $\tau'$ is an I-terminal edge, 
the sibling edges of $\tau'$ are O-edges 
which intersect $L_I$ by 
Lemma~\ref{LemSibling}(a). 
This contradicts the fact that 
all the edges intersecting $L_I$ 
are I-edges. 
Hence $L\subset L_I$.
Similarly 
we can show that 
if the label of $\tau$ is $s_I$, 
then $L\subset L_O$.

{\bf Statement (b)}.
Suppose that 
$\Gamma\cap ({\rm Int}~\Delta(\tau)-\tau)
\neq\emptyset$. 

Then
\begin{enumerate}
\item[(1)] 
Int~$\Delta(\tau)$  
contains a white vertex.
\end{enumerate}
Suppose that the label of $\tau$ is $s_O$. 
Let $N$ be a regular neighbourhood of 
$(e^*\cup \tau\cup e^{**})\cap D$ in $D$. 

Now $L\subset L_I$ by Statement (a).
Let $E=Cl(\Delta(\tau)-N)$ 
(see Fig.~\ref{fig-17}(b)). 
Since $L$ intersects only edges of label $s_I$, 
the tangle $(\Gamma\cap E,E)$ is 
an NS-tangle of label $s_I$ by (1). 
This contradicts Lemma~\ref{LemNS-Tangle}. 

Hence the label of $\tau$ is not $s_O$.
Similarly we can show that 
the label of $\tau$ is not $s_I$. 
Therefore we have 
$\Gamma\cap ({\rm Int}~\Delta(\tau)-\tau)
=\emptyset$.
\end{Proof}

Let $\Gamma$ be a minimal chart, and 
$(\Gamma\cap D,D)$ an N-tangle
with an IO-label pair $(s_I,s_O)$.
Let
$\mathbb T_O(D)=
\{\tau_1,\tau_2,\cdots,\tau_s\}$,
and 
$\mathbb T_I(D)=
\{\widetilde\tau_1,\widetilde\tau_2,
\cdots,\widetilde\tau_t\}$.
Let
$w_1,w_2,\cdots,w_s $ 
be the white vertices of
$\tau_1,\tau_2,\cdots,\tau_s$ 
respectively, 
and
$\widetilde w_1,\widetilde w_2,
\cdots,\widetilde w_t $ 
the white vertices of
$\widetilde\tau_1,\widetilde\tau_2,
\cdots,\widetilde\tau_t$ 
respectively.
Set
\begin{enumerate}
\item[]
$D^\dagger=
Cl(D-\cup_{i=1}^s\Delta(\tau_i)
-\cup_{j=1}^t\Delta(\widetilde\tau_j))$.
\end{enumerate}
The set $D^\dagger$ is called a 
{\it fundamental region} of the N-tangle.
By Lemma~\ref{lemDelta}, 
renumbering the terminal edges 
we always assume that
\begin{enumerate}
\item[(I)] 
the vertices
$\widetilde w_1,\widetilde w_2,\cdots,
\widetilde w_t $ are situated 
{\it clockwise} on 
$\partial D^\dagger$ in this order, and
\item[(II)] 
the vertices
$w_1,w_2,\cdots,w_s $ are situated 
{\it counterclockwise} on 
$\partial D^\dagger$ in this order.
\end{enumerate}
The tuple 
$(D^\dagger;
w_1,w_2,\cdots,w_s;
\widetilde w_1,\widetilde w_2,\cdots,
\widetilde w_t
)$
is called the {\it fundamental information}\index{fundamental information} 
for the N-tangle.

\begin{lem}\label{KeyLemma}
Let $\Gamma$ be a minimal chart, and 
$(\Gamma\cap D,D)$ an N-tangle with
an IO-label pair $(s_I,s_O)$.
Then the upward principal paths associated 
to the terminal edges of label $s_O$ in $D$ 
are mutually disjoint.
\end{lem}

\begin{figure}
\begin{center}
\includegraphics{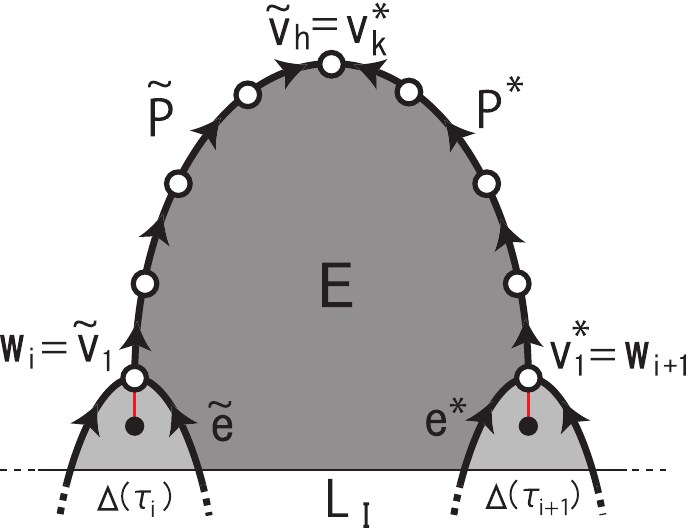}
\end{center}
\caption{ \label{fig-18}
The thick edges are of label $s_I$, and
the thin edges are of label $s_O$.}
\end{figure}

\begin{Proof}
Let
$(D^\dagger;
w_1,w_2,\cdots,w_s; 
\widetilde w_1,\widetilde w_2,
\cdots,\widetilde w_t)$ 
be the fundamental information 
for the N-tangle, and
$\mathbb T_O(D)=
\{\tau_1,\tau_2,\cdots,\tau_s\}$.

Now each upward principal path 
associated to a terminal edge 
of label $s_O$ in $D$ 
splits the set $D^\dagger$.
Thus to prove Lemma~\ref{KeyLemma},
it is sufficient to show
that for each $i=1,2,\cdots,s-1$,
two upward principal paths associated 
to the terminal edges $\tau_i$ 
and $\tau_{i+1}$  
are disjoint.

Suppose that for 
some integer $i$ $(1\le i<s)$,
there exist the 
two upward principal paths 
$\widetilde P,P^*$ associated to
the terminal edges 
$\tau_i,\tau_{i+1}$ 
with $\widetilde P\cap P^*\neq\emptyset$.

Let 
$(\widetilde v_1,\widetilde v_2,
\cdots,\widetilde v_p),
(v^*_1,v^*_2,\cdots,v^*_q)$
be vertex sequences of 
$\widetilde P,P^*$
respectively, 
here $\widetilde v_1= w_i,
v^*_1= w_{i+1}$. 
Let 
$h=\min\{ i~|~\widetilde v_i\in P^*\}$. 
Then $\widetilde v_{h}=v^*_{k}$ 
for some integer $0<k\le q$. 
Now $h\ge 2$, $k\ge 2$,
and 
$\widetilde P[\widetilde v_1,\widetilde v_{h}]
\cap P^*[v^*_1,v^*_{k}]=v^*_{k}$.
Let $(L_I,L_O)$ be a boundary IO-arc pair
of the N-tangle $(\Gamma\cap D,D)$,
and
$E$ the closure of 
the connected component of 
$D^\dagger-
(\widetilde P[\widetilde v_1,\widetilde v_{h}]
\cup P^*[v^*_1,v^*_{k}])$
with $E\cap\partial D\subset L_I$
(see Fig.~\ref{fig-18}).
Let $\widetilde e,e^*$
be the sibling edges of 
$\tau_i$ and $\tau_{i+1}$
intersecting $E-\{v_1^*,\widetilde v_1\}$ respectively. 
Set $\widetilde v_0$ and $v^*_0$ 
the vertices of 
$\widetilde e$ 
and 
$e^*$ outside $D$ respectively.
Then 
Remark~\ref{remarkPrincipal}(2) 
implies that 
the path with the vertex sequence 
$(\widetilde v_0,\widetilde v_1,
\cdots,\widetilde v_{h})$ is 
an upward-left-selective directed path,
and that
the path with the vertex sequence 
$(v^*_0,v^*_1,\cdots,v^*_{k})$ is 
an upward-right-selective directed path.
Thus $E$ is a half spindle by (II).
This contradicts Lemma~\ref{lemHalfSpindle}.
Therefore Lemma~\ref{KeyLemma} holds.
\end{Proof}

\begin{lem}\label{LemBijection}
Let $\Gamma$ be a minimal chart, and
$(\Gamma\cap D,D)$ an $N$-tangle 
with an IO-label pair $(s_I,s_O)$.
Then the map $f_D:\mathbb T_O(D)\to\mathbb T_I(D)$ 
is bijective.
\end{lem}

\begin{Proof}
Since the map $f_D:\mathbb T_O(D)\to\mathbb T_I(D)$
is injetive by Lemma~\ref{KeyLemma},
we have
\begin{enumerate}
\item[$(1)$]
$|\mathbb T_O(D)|\le |\mathbb T_I(D)|$.
\end{enumerate}

By changing all the orientations of 
the edges of the chart and
setting $\widetilde s_I=s_O,\widetilde s_O=s_I$, 
we obtain a new chart $\widetilde\Gamma$ and 
a new N-tangle 
$(\widetilde\Gamma\cap\widetilde D,\widetilde D)$
with
the IO-label pair $(\widetilde s_I,\widetilde s_O)$,
here $\widetilde D=D$.
Now
\begin{enumerate}
\item[$(2)$] 
all the terminal edges of label $s_I$ 
of the old N-tangle 
change to 
all the terminal edges 
of label $\widetilde s_O$ 
of the new N-tangle,
and
\item[$(3)$] 
all the terminal edges 
of label $s_O$ 
of the old N-tangle 
change to 
all the terminal edges 
of label $\widetilde s_I$ 
of the new N-tangle.
\end{enumerate}
Then for the new N-tangle, 
instead of $(1)$, 
we have 
$|\mathbb T_O(\widetilde D)|\le 
|\mathbb T_I(\widetilde D)|$.
Namely 
\begin{enumerate}
\item[$(4)$] 
$|\mathbb T_I(D)|\le |\mathbb T_O(D)|$.
\end{enumerate}
Thus we have 
$|\mathbb T_I(D)|=|\mathbb T_O(D)|$ 
by $(1)$ and $(4)$.
Therefore the map $f_D$ is bijective.
This proves Lemma~\ref{LemBijection}.\end{Proof}

{\bf Proof of Theorem~\ref{Theorem1}.}
Let $\Gamma$ be a minimal chart, and 
$(\Gamma\cap D,D)$ a net-tangle with 
a label pair $(\alpha,\alpha+1)$.
Then the tangle is an N-tangle by
Lemma~\ref{LemChromaticTangleIsN-tangle}(b).
Thus Theorem~\ref{Theorem1}(a) holds. 
Theorem~\ref{Theorem1}(c) follows from 
Lemma~\ref{LemChromaticTangleIsN-tangle}(a).
Also Theorem~\ref{Theorem1}(d) holds by 
Lemma~\ref{LemBijection}.
Thus $T_I(D)=T_O(D)$ by 
Lemma~\ref{LemTerminalEdgeInN-Tangle}. 
Hence we have $E_I(D)=E_O(D)$ by 
Lemma~\ref{EI+TI=EO+TO}.
Thus Theorem~\ref{Theorem1}(b) holds.
Finally, Theorem~\ref{Theorem1}(e) holds
by Lemma~\ref{LemTerminalEdgeInN-Tangle}.
This proves Theorem~\ref{Theorem1}.{\hfill {$\square$}}

\section{{\large Proof of Theorem~1.2}}
\label{s:ProofTh2}

In this section
we prove Theorem~\ref{Theorem2}.

Let $\Gamma$ be a chart.
For a subset $X$ of $\Gamma$,
let
\begin{enumerate}
\item[]
$B(X)=$ the union of all the disk 
bounded by a cycle in $X$, and
\item[]
$T(X)=$ the union of all the terminal edge 
intersecting $X\cup B(X)$.
\end{enumerate}
The set $X\cup B(X)\cup T(X)$ is called
the {\it SC-closure}\index{SC-closure} of $X$ and denoted by $SC(X)$.

\begin{remark}\label{RemSC-Closure}
{\rm Each connected component of 
the SC-closure $SC(X)$ is simply connected.}
\end{remark}

\begin{figure}
\begin{center}
\includegraphics{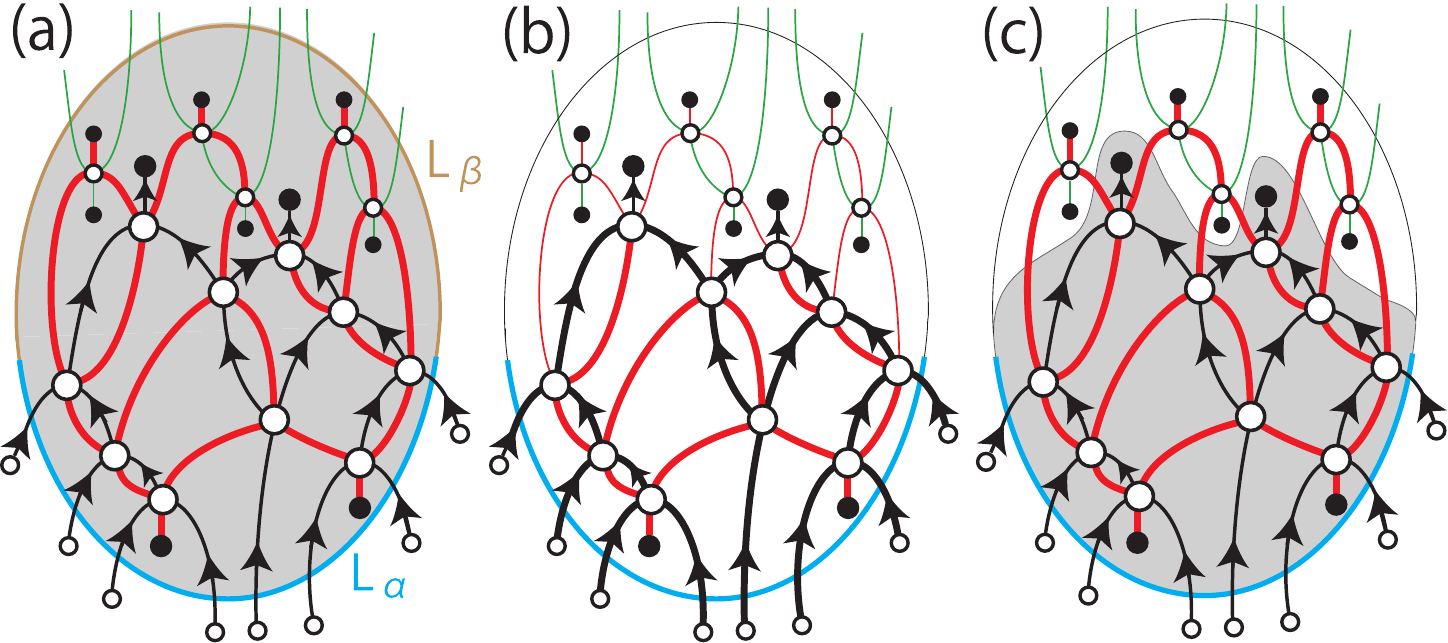}
\end{center}
\caption{ \label{fig-19}
(a) the gray area is the disk $D$, 
the edges with arrows are of label $\alpha$,
the thick edges are of label $\alpha+1$,
the thin edges are of label $\alpha+2$.
(b) the big white vertices are in $\mathcal V$,
the union of the thick edges is 
$\mathcal E_1\cup \mathcal E_2$.
(c) the gray area is the disk $N(\Gamma,D,\alpha)$.}
\end{figure}

Let $\Gamma$ be a chart,
and 
$(\Gamma\cap D,D)$ a net-tangle
with a label pair $(\alpha,\beta)$ and 
a boundary arc pair $(L_\alpha,L_\beta)$.
Let (see Fig.~\ref{fig-19})
\begin{enumerate}
\item[]
$\mathcal V =\cup
\{v~|~\text{$v$ is a vertex in 
$D\cap\Gamma_\alpha\cap\Gamma_{\alpha+1}$}
\}$,
\item[]
$\mathcal E_1=\cup
\{e~|~
\text{$e$ is an edge of label $\alpha$
intersecting $L_\alpha$}
\}$,
\item[]
$\mathcal E_2=\cup
\{e~|~
\text{$e$ is an edge connecting 
two white vertices in $\mathcal V$}
\}$,
\item[]
$Y=$the connected component of
$(L_\alpha\cup\mathcal V\cup
\mathcal E_1\cup\mathcal E_2)\cap D$
containing $L_\alpha$,
\item[]
$N(\Gamma,D,\alpha)=$a regular neighbourhood 
of the SC-closure $SC(Y)$ in $D$,
\item[]
$L^*_{\alpha}=
Cl(\partial D\cap N(\Gamma,D,\alpha))$,
and
\item[]
$L^*_{\alpha+1}=
Cl(\partial N(\Gamma,D,\alpha)-
L^*_{\alpha})$.
\end{enumerate}
Then $N(\Gamma,D,\alpha)$ is simply connected
by Remark~\ref{RemSC-Closure}.
Thus the set $N(\Gamma,D,\alpha)$ is a disk.
The number  
$|\Gamma_{\alpha+1}\cap 
\partial N(\Gamma,D,\alpha)|$
is called the $(D,\alpha+1)$-{\it number} \index{$(D,\alpha+1)$-number $n(\Gamma,D,\alpha+1)$}
of the chart 
and denoted by $n(\Gamma,D,\alpha+1)$.
Let 
\begin{enumerate}
\item[]
$N(L^*_{\alpha+1})=$  
a regular neighbourhood of $L^*_{\alpha+1}$ in $D$, and
\item[]
 $\mathcal{S}(\Gamma,D,\alpha+1)=$ 
the set of all 
the chart obtained from $\Gamma$ 
by a finite sequence of 
C-I-M2 moves in $N(L^*_{\alpha+1})$.
\end{enumerate}
A chart in $\mathcal S(\Gamma,D,\alpha+1)$ is 
$(\Gamma,D,\alpha+1)$-{\it minimal}
if its $(D,\alpha+1)$-number is minimal\index{$(\Gamma,D,\alpha+1)$-minimal}
amongst $\mathcal{S}(\Gamma,D,\alpha+1)$. 

\begin{remark}\label{RemSCN}
{\rm 
Let $\Gamma$ be a chart,
and 
$(\Gamma\cap D,D)$ a net-tangle
with a label pair $(\alpha,\beta)$ and 
a boundary arc pair $(L_\alpha,L_\beta)$.
\begin{enumerate}
\item[$(1)$] 
If an edge intersects 
$\partial N(\Gamma,D,\alpha)$,
then the edge transversely intersects 
$\partial N(\Gamma,D,\alpha)$.
\item[$(2)$] 
$L^*_{\alpha+1}$ is 
a proper arc of $D$.
\item[$(3)$] $\partial N(\Gamma,D,\alpha)=
L^*_\alpha\cup L^*_{\alpha+1},$ 
$\Gamma\cap L^*_\alpha=
\Gamma_\alpha\cap L^*_\alpha,$ and 
$\Gamma\cap L^*_{\alpha+1}=
\Gamma_{\alpha+1}\cap L^*_{\alpha+1}$.
\item[$(4)$] 
Any chart in $\mathcal S(\Gamma,D,\alpha+1)$ 
is M2-related to $\Gamma$.
\end{enumerate}
}
\end{remark}

{\bf Proof of Theorem~\ref{Theorem2}}.
We prove Theorem~\ref{Theorem2}
by induction on the number $|\alpha-\beta|$.
Let $\Gamma$ be a minimal chart, and
$(\Gamma\cap D,D)$ a net-tangle with 
a label pair $(\alpha,\beta)$ and
a boundary arc pair $(L_\alpha,L_\beta)$.

Suppose $|\alpha-\beta|=1$.
Then the tangle is a net-tangle
with a label pair $(\alpha,\alpha+1)$.
Thus it is an N-tangle by
Lemma~\ref{LemChromaticTangleIsN-tangle}(b).
Hence Theorem~\ref{Theorem2} holds 
for the case $|\alpha-\beta|=1$.

Suppose that $|\alpha-\beta|>1$.
Without loss of generality we can assume that
\begin{enumerate}
\item[$(1)$]
all the edges of label $\alpha$ 
intersecting $\partial D$
are I-edges for $D$.
\end{enumerate}
We use all the notations in the definition of 
$(\Gamma,D,\alpha+1)$-minimal. 
Now 
$\Gamma\cap\partial N(\Gamma,D,\alpha)=
(\Gamma\cap L^*_\alpha)\cup
(\Gamma\cap L^*_{\alpha+1})
\subset
\Gamma_\alpha\cup\Gamma_{\alpha+1}$
by Remark~\ref{RemSCN}$(3)$. 
Thus by Boundary Condition Lemma 
(Lemma~\ref{BoundaryConditionLemma})
we have
\begin{enumerate}
\item[$(2)$]
$\Gamma\cap N(\Gamma,D,\alpha)
\subset
\Gamma_\alpha\cup\Gamma_{\alpha+1}$.
\end{enumerate}

{\bf Claim~1}.
Any edge of label $\alpha+1$ intersecting 
$\partial N(\Gamma,D,\alpha)$ 
is not an edge of a bigon.

{\it Proof of Claim~$1$}.
Suppose that 
an edge $e$ of label $\alpha+1$ intersecting 
$\partial N(\Gamma,D,\alpha)$ is 
an edge of a bigon.
Then
\begin{enumerate}
\item[$(3)$] 
the edge $e$ transversely intersects
$\partial N(\Gamma,D,\alpha)$ by Remark~\ref{RemSCN}(1).
\end{enumerate}
Since $|\alpha-\beta|>1$, we have 
$\alpha<\alpha+1<\beta$. 
Thus
\begin{enumerate}
\item[$(4)$] $e\subset D$.
\end{enumerate}
Let $e^*$ be the other edge of the bigon.
Then $e^*$ intersects $N(\Gamma,D,\alpha)$.
Further $e^*$ is of label $\alpha$ or $\alpha+2$,
because 
$e$ is of label $\alpha+1$.
Since no edge of label $\alpha+2$ intersects
$N(\Gamma,D,\alpha)$ by $(2)$,
the edge $e^*$ is of label $\alpha$.
Hence
the two white vertices of $e$ 
are in 
$\Gamma_\alpha\cap\Gamma_{\alpha+1}$.
Hence $e\subset \mathcal E_2
\subset N(\Gamma,D,\alpha)$ by $(4)$.
This contradicts Statement $(3)$. 
Hence Claim~$1$ holds.

By Claim 1 all the chart in 
$\mathcal S(\Gamma,D,\alpha+1)$ 
is a minimal chart. 
Without loss of generality we can assume that
\begin{enumerate}
\item[$(5)$]
the chart $\Gamma$ is 
$(\Gamma,D,\alpha+1)$-minimal.
\end{enumerate}

{\bf Claim~2}. 
All the edges intersecting $L^*_{\alpha+1}$ 
are I-edges for $N(\Gamma,D,\alpha)$ 
of label $\alpha+1$ or
all the edges intersecting $L^*_{\alpha+1}$ 
are O-edges for $N(\Gamma,D,\alpha)$
of label $\alpha+1$.

{\it Proof of Claim~$2$}.
Suppose that there exist two edges $e_1,e_2$ 
of label $\alpha+1$ intersecting $L^*_{\alpha+1}$
such that
one of $e_1,e_2$ is an I-edge 
for $N(\Gamma,D,\alpha)$
and 
the other an O-edge for $N(\Gamma,D,\alpha)$.
Let $L$ be the subarc of $L^*_{\alpha+1}$
connecting the two edges $e_1$ and $e_2$. 
Without loss of generality
we can assume that 
$\Gamma_{\alpha+1}\cap{\rm Int}~L=\emptyset$.
Let $\widetilde\Gamma$ be a chart 
obtained from $\Gamma$
by applying a C-I-M2 move between 
$e_1$ and $e_2$ along the arc $L$ 
in $N(L^*_{\alpha+1})$.
Since each of $e_1$ and $e_2$ contains 
two white vertices,
$\widetilde\Gamma\cap D$ does not contain
a free edge nor a hoop.
Further Claim~$1$ assures that 
$\widetilde\Gamma$ is a minimal chart 
with
$n(\widetilde\Gamma,D,\alpha+1)<
n(\Gamma,D,\alpha+1)$.
This contradicts Statement $(5)$.
Thus Claim~$2$ holds.

By Claim~$2$, the tangle
$(\Gamma\cap N(\Gamma,D,\alpha),N(\Gamma,D,\alpha))$ 
is a net-tangle with 
$\Gamma\cap N(\Gamma,D,\alpha)
\subset\Gamma_\alpha\cup\Gamma_{\alpha+1}$.
Thus it is an N-tangle by 
Lemma~\ref{LemChromaticTangleIsN-tangle}(b).
Hence it is upward by $(1)$.
Let $D^*=Cl(D-N(\Gamma,D,\alpha))$ and 
$L^*_\beta=Cl(\partial D^*-L^*_{\alpha+1})$.
Then $\Gamma\cap\partial D^*\subset
\Gamma_{\alpha+1}\cup\Gamma_\beta$.
Thus $\Gamma\cap D^*\subset
\cup_{i=\alpha+1}^\beta\Gamma_i$
by Boundary Condition Lemma 
(Lemma~\ref{BoundaryConditionLemma}).
Thus $(\Gamma\cap D^*,D^*)$ is 
a net-tangle with 
a label pair $(\alpha+1,\beta)$
and 
a boundary pair $(L^*_{\alpha+1},L^*_\beta)$.
Since $|(\alpha+1)-\beta|<|\alpha-\beta|$,
the results follow 
by the induction.
This proves Theorem~\ref{Theorem2}.{\hfill {$\square$}}

\newpage

{\bf List of terminologies}\vspace{2mm}\\
{\small $
\begin{array}{ll||}
\text{2-color disk} & p13 \\
\text{3-color disk} & p13 \\
\text{arc triplet} & p26 \\
\text{bigon} & p8 \\
\text{boundary arc pair} & p5 \\
\text{boundary IO-arc pair} & p20 \\
\text{C-move~equivalent} & p1 \\
\text{consecutive triplet} & p12 \\
\text{corresponding terminal edge for $\tau_O$} & p31 \\
\text{cycle} & p13 \\
\text{dichromatic path} & p14 \\
\text{directed cycle} & p15 \\
\text{directed path} & p14 \\
\text{dominate} & p17 \\
\text{downward-left-selective} & p17 \\
\text{downward maximal} & p17 \\
\text{downward net-tangle} & p5 \\
\text{downward-right-selective} & p17 \\
\text{edge sequence} & p14 \\
\text{free edge} & p8 \\
\text{fundamental information} & p33 \\
\text{half spindle} & p26 \\
\text{hoop} & p9 \\
\text{I-edge} & p3 \\
\text{inside edge} & p23 \\
\text{inward} & p3, p16 \\
\text{IO-label pair} & p20 \\
\text{IS-arc} & p24 \\
\text{I-terminal edge} & p5 \\
\text{label pair} & p5, p23, p26 \\
\text{leading to $e$} & p17 \\
\text{locally left-side} & p16 \\
\text{locally right-side} & p16 \\
\end{array}
$
~~
$\begin{array}{ll}
\text{loop} & p15 \\
\text{M$\&$M path} & p14 \\
\text{M2-related to a chart} & p6 \\
\text{middle arc} & p7 \\
\text{middle at $v$} & p7 \\
\text{minimal chart} & p9 \\
\text{net-tangle} & p4 \\
\text{NS-tangle} & p13 \\
\text{N-tangle} & p5 \\
\text{O-edge} & p3 \\
\text{OS-arc} & p24 \\
\text{O-terminal edge} & p5 \\
\text{outside edge} & p13, p23 \\
\text{outward} & p3, p16 \\
\text{oval nest} & p9 \\
\text{path of label $m$} & p14 \\
\text{path pair} & p23 \\
\text{principal path} & p 30, p31 \\
\text{proper arc} & p6 \\
\text{ring} & p9 \\
\text{SC-closure $SC(X)$} & p35 \\
\text{sibling edge} & p19 \\
\text{side-edge} & p14 \\
\text{simple hoop} & p9 \\
\text{spindle} & p22 \\
\text{starting form $e$} & p 17\\
\text{tangle} & p2 \\
\text{terminal edge} & p5 \\
\text{upward-left-selective} & p17 \\
\text{upward maximal} & p17 \\
\text{upward net-tangle} & p5 \\
\text{upward-right-selective} & p17 \\
\text{vertex sequence} & p14 \\
\end{array}
$}

\vspace{1cm}
{\bf List of notations}\vspace{2mm}\\
{\small $
\begin{array}{ll||}
\text{$\Gamma_m$} & p2 \\
\text{${\rm Main}(\Gamma)$} & p12 \\
\text{${\mathcal{W}}_O^{{\rm Mid}}(C,m)$} & p13 \\
\text{$P[v_i,v_j]$} & p23 \\
\text{$E_I(D)$, $E_O(D)$} & p29 \\
\end{array}
$
~~
$\begin{array}{ll}
\text{$T_I(D)$, $T_O(D)$} & p29 \\
\text{$\mathbb T_I(D)$, $\mathbb T_O(D)$} & p31 \\
\text{$(\Gamma,D,\alpha+1)$-minimal} & p35 \\
\text{$(D,\alpha+1)$-number $n(\Gamma,D,\alpha+1)$} & p35 \\
& \\
\end{array}
$}

\end{document}